\documentclass[reqno,11pt]{amsart}  

\usepackage{gensymb}
 
 \usepackage{amsmath}
\usepackage{amssymb} 
\usepackage{mathtools}   
\usepackage{graphicx}
\usepackage{graphics}  
\usepackage{color}
\usepackage{verbatim} 
\usepackage[normalem]{ulem} 
\usepackage[mathscr]{eucal}
\usepackage{upgreek}
\usepackage{enumerate} 
\usepackage[titletoc]{appendix}
\numberwithin{equation}{section}  
\usepackage[refpage,noprefix]{nomencl}
\usepackage{nomencl} 

\makenomenclature

\usepackage{bbm}

\usepackage[titletoc]{appendix}
\numberwithin{equation}{section}  
 
\usepackage[refpage,noprefix]{nomencl}
\usepackage{nomencl} 

\newtheorem{theorem}{Theorem}[section] 
\newtheorem{lemma}[theorem]{Lemma} 
\newtheorem{proposition}[theorem] {Proposition} 
 
\newtheorem{remark}[theorem]  {Remark} 
\newtheorem{definition}[theorem] {Definition}

\theoremstyle{definition}
\newtheorem{example}[theorem] {Example}

 \usepackage{float}
\restylefloat{figure}

\DeclareMathAlphabet{\mathpzc}{OT1}{pzc}{m}{it}

\DeclarePairedDelimiter{\abs}{\lvert}{\rvert}
\DeclarePairedDelimiter{\norm}{\lVert}{\rVert}

\renewcommand{\L} {\Lambda} %

\def\d{\delta}

\newfam\Bbbfam 
\font\tenBbb=msbm10 
\font\sevenBbb=msbm7 
\font\fiveBbb=msbm5 
\textfont\Bbbfam=\tenBbb 
\scriptfont\Bbbfam=\sevenBbb 
\scriptscriptfont\Bbbfam=\fiveBbb

 \newcommand{\C}     {\mathbb{C}} 
\newcommand{\R}     {\mathbb{R}} 
\newcommand{\Z}     {\mathbb{Z}} 
\newcommand{\N}     {\mathbb{N}} 
\renewcommand{\P}   {\mathbb{P}} 
 
\newcommand{\E}     {\mathbb{E}} 
 
\newcommand{\T}     {\mathbb{T}}

\def\1{{\mathchoice {1\mskip-4mu\mathrm l}      
{1\mskip-4mu\mathrm l} 
{1\mskip-4.5mu\mathrm l} {1\mskip-5mu\mathrm l}}} 
\newcommand{\ssup}[1] {{\scriptscriptstyle{({#1}})}} 
\def\comment#1{} 
\newtheoremstyle{thm}{2ex}{2ex}{\itshape\rmfamily}{} 
{\bfseries\rmfamily}{}{1.7ex}{} 
 
\newtheoremstyle{rem}{1.3ex}{1.3ex}{\rmfamily}{} 
{\itshape\rmfamily}{}{1.5ex}{} 
 
\newenvironment{proofsect}[1] 
{\vskip0.1cm\noindent{\scshape #1.}\hskip0.5cm}

\newcommand{\bG} {\boldsymbol{G}}

\newcommand{\bL} {\boldsymbol{L}}

\newcommand{\bV} {\boldsymbol{V}}

\newcommand{\bv} {\boldsymbol{v}}

\newcommand{\bpi}{\boldsymbol{\pi}}
\newcommand{\bomega}{\boldsymbol{\omega}}

\newcommand{\bLscr}{\boldsymbol{\Lscr}}
 \newcommand{\bGamma}{\boldsymbol{\Gamma}}
\newcommand{\gb}{{\mathbf{\mathcal{G}}}}
\newcommand{\Lb}{{\mathbf{L}}}

\newcommand{\bphi}{\boldsymbol{\phi}} 
 \newcommand{\bmu}{\boldsymbol{\mu}}

\newcommand{\Acal}   {{\mathcal A }}

\newcommand{\Fcal}   {{\mathcal F }} 
\newcommand{\Gcal}   {{\mathcal G }}

\newcommand{\Lcal}   {{\mathcal L }} 
\newcommand{\Mcal}   {{\mathcal M }} 
 
\newcommand{\Ocal}   {{\mathcal O }}

\newcommand{\Scal}   {{\mathcal S }}

\newcommand{\Hscr} {\mathscr{H}}

\newcommand{\Lscr} {\mathscr{L}}

\newcommand{\Ex}{\mathsf{Exp}}

\def\Poi{{\mathsf P}}

\renewcommand{\Pr} {\mathsf{P}}
\newcommand{\Er} {\mathsf{E}}
\newcommand{\Vsf}{\mathsf{V}}
\newcommand{\Jsf}{\mathsf{J}}
\newcommand{\psf}{\mathsf{p}}

 \newcommand{\ex}{{\rm e}} 
 
\renewcommand{\d}{{\rm d}}

\newcommand{\Sym}{\mathfrak{S}}
\newcommand{\id}{{\sf{I}}}

\newcommand{\supp}{{\operatorname {supp}}}

\newcommand{\diag}{{\operatorname {diag}\,}}

\newcommand{\wind}{{\operatorname {\sf wind}}}
\newcommand{\tr}{{\operatorname {Tr}}}

\newcommand{\Exp}{\mathscr{E}\kern-0.2mm{\operatorname{xp}}}
\newcommand{\Log}{\mathscr{L}\kern-0.2mm{\operatorname{og}}}

\newcommand{\heap}[2]{\genfrac{}{}{0pt}{}{#1}{#2}} 
\renewcommand{\emptyset} {\varnothing}

\begin{document}

\title[\hfill Bosonic loop measures \hfill]
{Space-time random walk loop measures} 
 

\author{Stefan Adams  and  Quirin  Vogel}
\address{Mathematics Institute, University of Warwick, Coventry CV4 7AL, United Kingdom}
\email{S.Adams@warwick.ac.uk, Q.Vogel@warwick.ac.uk}

\thanks{}
  

 


%
%

\begin{abstract} 
In this  work, we introduce and investigate two novel classes of loop measures, space-time Markovian loop measures and  \textit{Bosonic loop measures}, respectively. We consider  loop soups with intensity $ \mu\le 0 $ (chemical potential in physics terms), and secondly, we study Markovian loop measures on graphs with an additional ``time'' dimension leading to so-called space-time random walks and their loop measures and  Poisson point loop processes. Interesting phenomena appear when the additional coordinate of the space-time process is on a discrete torus with non-symmetric jump rates. The projection of these space-time random walk loop measures onto the space dimensions is loop measures on the spatial graph,  and in the scaling limit of the discrete torus, these loop measures converge to the so-called \textit{Bosonic loop measures}.  This  provides a natural probabilistic definition of  \textit{Bosonic loop measures}. These novel loop measures have similarities with the standard Markovian loop measures only that they  give weights to loops of certain lengths, namely any length which is multiple of a given length $ \beta> 0 $ which serves as an additional parameter. 
We complement our study with generalised versions of Dynkin's isomorphism theorem (including a version for the whole complex field) as well as Symanzik's moment formulae for complex Gaussian measures. Due to the lacking symmetry of our space-time random walks, the distributions of the occupation time fields are  given in terms of  complex Gaussian measures over complex-valued random fields (\cite{B92,BIS09}). Our space-time setting allows obtaining quantum correlation functions as  torus limits of space-time correlation functions.

  \end{abstract} 
 \maketitle 
\section{Introduction}\label{S:intro}
Recently Markovian Loop measures have become an active field in probability theory with its origin going back to Symanzik \cite{Sym69}.  Our focus in this work is  Markovian loop measures on graphs with an additional ``time'' dimension leading to so-called space-time random walks and their loop measures and  Poisson point loop processes.  This work will show that interesting new phenomena appear when the added coordinate of the space-time process is on a discrete torus with non-symmetric jump rates. The projection of these space-time random walk loop measures onto the space dimensions is loop measures on the spatial graph,  and in the scaling limit of the discrete torus, these loop measures converge to the  Bosonic loop measures which are a new class of loop measures. These have different properties than the Markovian loop measures some of which we study in this work. The Bosonic loop measures not only have the probabilistic derivation as torus limits of space-time Markovian loop measures, a second major interest in these objects stems from the fact that the total weight of the Bosonic loop measure for a finite graph is exactly the logarithm of the grand-canonical partition function of  non-interacting Bose gas on a finite graph  in  thermodynamic equilibrium at an inverse temperature $ \beta>0 $ and chemical potential $ \mu\le 0 $.  The study of Markovian loop measures has been outlined in the lecture notes \cite{LeJan} and \cite{Sznitman} with more  recent developments  in \cite{FR14}.

\subsection{Notation and set up}
We begin with the definition of loop measures for processes with discrete state space. Let $\gb$ be a finite set. We will later specify our choice of $ \gb$ (usually endowed with a graph structure). Let $W=(w(x,y))_{x,y\in\gb}$ be an irreducible matrix (not necessarily symmetric here) with non-negative entries indexed by $\gb$. We say that $x,y\in\gb$ are connected (denoted by $x\sim y$) if $w(x,y)>0$, and we set $w(x,x)=0 $. Let $\kappa=(\kappa(x))_{x\in\gb}$ be a vector with non-negative entries. We often refer to $\kappa$ as \textsf{killing} if $\kappa(x)>0$ for some $x\in\gb $. Denote the normalising constant $\lambda(x)=\kappa(x)+\sum_{y\sim x}w(x,y)$. Clearly $P=(p({x,y}))_{x,y\in\gb}$ defined by $p(x,y)=\tfrac{w(x,y)}{\lambda(x)}$ induces (sub-) Markovian transition probabilities on $\gb$ via the generator matrix $Q=(q(x,y))_{x,y\in\gb}$ given by $ Q=\lambda(P-I) $. Here, $\lambda=\texttt{diag}(\lambda)$ is the matrix with $(\lambda(x))_{x\in\gb}$ on the diagonal and all other entries set to zero. For $t\ge 0$ transition densities are given by 
\begin{equation}
p_t(x,y)=\ex^{tQ}(x,y)\quad\, x,y\in\gb.
\end{equation}
The transition densities define a continuous time Markov process on the Skorohod space of c\`adl\`ag paths from $[0,\infty)\to\gb$ when adding a cemetery state (necessary if $\kappa\neq 0$). For $x\in\gb$, $\P_x$ denotes the distribution of that (sub-)Markov process with initial condition $x$, i.e., $\P_x(X_0=x)=1$.  We call $ X=(X_t)_{t\ge 0} $ the random walk process on the graph $\gb $ with killing $ \kappa$. When $ \gb$ is a finite connected subset of $\Z^d$ we write  $\gb=\L\Subset\Z^d $, and we will sometimes consider killing upon reaching the boundary, i.e., $ \kappa(x)=\sum_{y\in\Z^d\setminus\L}\1\{x\sim y\}, x\in\L$ in the case of unit weights, which corresponds to Dirichlet boundary conditions.

The random walk process defines different measures on the space of loops. Following \cite{Sznitman}, for any $ t>0 $ let $ \Gamma_t $ be the space of c\`adl\`ag functions $ \omega\colon[0,t]\to\gb $ with the same value at $0$ and $t$.  We denote by $X_s, 0\le  s\le t $, the canonical coordinates, and we extend them to be equal to a "graveyard" state for $s>t$, so that $X_s(\omega) $  is well-defined for any $ s\in\R$. Note that the spaces  $\Gamma_t $ are pairwise disjoint, as $t$ varies over $(0,\infty) $, and define
 \begin{equation}
\Gamma:=\bigcup_{t>0} \;\Gamma_t\, .
\end{equation}
For each $ \omega\in\Gamma $ we denote the time horizon or the length of the loop $\omega$ by $ \ell(\omega) $ being the unique $ t>0 $ such that $ \omega\in\Gamma_t $. The $ \sigma$-algebra $\bG$ is obtained from the bijection $ \Gamma_1\times (0,\infty) \to\Gamma $ defined as
$$
\Gamma_1\times(0,\infty)\ni(\omega,t)\mapsto \omega(\cdot)=\omega\big(\frac{\cdot}{t}\big)\in\Gamma,
$$
where we endow $ \Gamma_1\times(0,\infty) $ with the canonical product $\sigma$-algebra (and $ \Gamma_1 $ is endowed with the $ \sigma$-algebra generated by $X_t, 0\le t\le 1$). The $ \sigma $-algebra $\bG$ on $ \Gamma $ is the image of the $ \sigma$-algebra on $ \Gamma_1\times(0,\infty) $. We define for $ x,y\in\gb, t>0$, the measure $ \P^{\ssup{t}}_{x,y} $ as the image of $\1\{X_t=y\}\P_x$ under the map $ (X_s)_{0\le s\le t} $ from $ D_\gb\cap\{X_t=y\} $ into itself; if $x=y$ we identify it as a measure on $ \Gamma_t $. Here $ D_{\gb} $ is the Skorohod space over $ \gb$ with time horizon $ [0,\infty) $.  We study the following new loop measures which are related  to the well-known Markovian loop measures (see e.g. \cite{LeJan} for an overview on Markovian loop measures and \cite{CS16} for introducing a chemical potential).

\begin{definition}[\textbf{Loop measures}]
Let $C\subset \Gamma$ be measurable, $\mu\le 0$, and $\beta>0$. The \textit{Markovian loop measure} $M_{\gb,\mu}$ on $ \gb $ with chemical potential $ \mu $ is defined as
\begin{equation}
M_{\gb,\mu}[C]=\sum_{x\in\gb}\int_0^\infty\frac{\ex^{t\mu}}{t}\P_{x,x}^{\ssup{t}}\left(C\right)\d t=\sum_{x\in\gb}\int_0^\infty\frac{\ex^{t\beta\mu}}{t}\P_{x,x}^{\ssup{\beta t}}\left(C\right)\d t.
\end{equation}
The  \textit{Bosonic loop measure} $ M_{\gb,\mu,\beta}^B $ on $ \gb $ with chemical potential $ \mu $ and time horizon $ \beta $ is defined as 
\begin{equation}\label{BosLoopMDef}
M_{\gb,\mu,\beta}^B[C]=\sum_{x\in\gb}\sum_{j=1}^\infty\frac{\ex^{\beta j\mu}}{j}\P_{x,x}^{\ssup{j\beta}}\left(C\right)\,.
\end{equation}
Note that both measures are sigma-finite since the mass of $ \big(\bigcup_{t=1/r}^r\Gamma_t\big) $ is finite for any $ r\in\N $ under both measures. 
\end{definition}

The Markovian loop measure $ M_{\gb,\mu} $ with chemical potential  $ \mu\le 0 $ puts lower weight on loops with larger time horizons. The main novelty are the Bosonic loop measures whose probabilistic meaning is revealed in Section~\ref{Sec-2} as limits of space-time random walk loop measures (Theorem~\ref{THM-main1}). The major significances of these loop measures are explained in Section~\ref{sec-bosons}, in particular one can show that the logarithm of the mass of the Bosonic loop measure is just the so-called pressure of some system of Bosons, see  \eqref{partition3}.

\subsection{Overview}
The main focus in this work is the Bosonic loop measure $ M_{\gb,\mu,\beta}^B $ which is defined in similar way to the Markovian loop measure (\cite{LeJan}) but differs in having support only on loops with time horizons being a multiple of the given parameter $ \beta $.  The Bosonic loop measures for simple random walks are linked to  equilibrium quantum statistical mechanics, see Section~\ref{S:Bosons} for details.  The main idea and novelty of our study is to view these measures as natural Markovian loop measures in a space-time structure where we add an additional ``time" torus to the given graph.  After providing technical tools and different classes of space-time random walks in Section~\ref{Markovian and Bosonic loop measures},  we prove the first main result about the torus limit approximation of the Bosonic loop measure in Theorem~\ref{THM-main1}. In Theorem~\ref{THM-main2} we deliver the corresponding torus limit for the distribution of the occupation time field of the underlying Poisson loop processes, where we in particular distinguish between real torus jumps and solely spatial jumps. We finish Section~\ref{Markovian and Bosonic loop measures} with  various examples for our three different classes of space-time random walks, and we finally show in Theorem~\ref{support_prop} that certain space-time random walk loop measures have support which equal the ones of the Bosonic loop measures.

Our main results concern isomorphism theorems (Theorem~\ref{Isomorphism occupation field}) for the square of the complex field as well as the complex field and generalised Symnazik's type formulae for moments of the occupation time fields, see  Proposition~\ref{Symanzik_formula_generalized} and Theorem~\ref{Space-time loop-representation of the correlation function}. Taking the so-called torus limit, Theorem~\ref{Space-time loop-representation of the correlation function} allows to represent quantum correlation functions as space-time Green functions of the corresponding complex Gaussian measure. Thus we deliver a purely probabilistic derivation of Bosonic loop measures and quantum correlation functions.  We propose to analyse these complex measures in greater detail in the future to study Bose-Einstein condensation phenomena and their possible connection to random interlacements (\cite{Sznitman}) and to general permanental processes (\cite{FR14}).  In Section~\ref{sec-bosons} we give a review about  Boson systems and prove  Poisson loop process representations of the partition function and the quantum correlation functions.  The Bosonic loop measure and its Poisson loop process is a natural extension of random walk loop soups (\cite{LP14,LF07}), and it seems feasible to define Bosonic loop measures and soups also over the continuum space $\R^d $, e.g., see \cite{LW04}. In Section~\ref{S-limit} we briefly outline which of our results  hold in the so-called thermodynamic limit $ \L\uparrow\Z^d $ when $\gb=\L\Subset\Z^d$.

\section{Results}\label{Sec-2}
We state properties of our loop measures in a finite graph setting in Section~\ref{S-properties} and then show that the space-time Markovian loop measures converge to the Bosonic ones in the torus limit. In Section~\ref{Isomorphism Theorems and space-time correlation functions}, we study  so-called isomorphism theorems for finite graphs and extend them to  non-symmetric settings which arise naturally for our space-time loop measures. This allows us to obtain moments formulae for perturbed complex Gaussian measures as well as space-time loop representations of quantum correlation functions. In Section~\ref{S-limit} we briefly discuss  taking limits towards countable infinite graphs.
\subsection{Markovian and Bosonic loop measures}\label{Markovian and Bosonic loop measures}\label{S-properties}

We collect some basic results about the Markovian loop measures  and the Bosonic loop measure.

\begin{proposition}[\textbf{Properties}]\label{properties}
Let $ \mu\le 0 $ and $  \beta>0 $. For the Bosonic loop measure the  finite-dimensional distributions for $ 0< t_1<t_2<\cdots<t_k <\infty, k\in\N $, and $A\subset(0,\infty)$ measurable with $ \inf\{A\}>t_k $ are
\begin{equation}\label{finite_Boson _distr}
\begin{split}
&M_{\gb,\mu,\beta}^B\left[X_{t_1}=x_1,\ldots,X_{t_k}=x_k,\,\ell\in A\right]=\\[1.5ex]
&\sum_{ j=1}^\infty p_{t_2-t_1}(x_1,x_2)\cdots p_{t_k-t_{k-1}}(x_{k-1},x_k)p_{\beta j -t_k+t_1}(x_k,x_1)\frac{\ex^{\beta\mu j}}{j}\mathbbm{1}_A(\beta j)\, \,\mbox{ when } k>1,
\end{split}
\end{equation} 
and
\begin{equation}\label{finite_Boson _distr k=1}
\begin{split}
M_{\gb,\mu,\beta}^B\left[X_{t_1}=x_1,\,\ell\in A\right]=\sum_{j=1}^\infty p_{\beta j}(x_1,x_1)\frac{\ex^{\beta\mu j}}{j}\mathbbm{1}_A(\beta j)\, \mbox{ when } k=1.
\end{split}
\end{equation} 
The finite-dimensional distributions for the Markovian loop measures are similar by just  replacing the sum in equation \eqref{finite_Boson _distr} and \eqref{finite_Boson _distr k=1} by the integral such that the support of the Markovian loop measure is $ \Gamma $ instead of the union $ \bigcup_{j=1}^\infty \Gamma_{j\beta} $,
\begin{equation}\label{FDD-boson}
\begin{split}
&M_{\gb,\mu}\left[X_{t_1}=x_1,\ldots,X_{t_k}=x_k,\,\ell\in A\right]=\\[1ex]
&\int_{A} p_{t_2-t_1}(x_1,x_2)\cdots p_{ t_k-t_{k-1}}(x_{k-1},x_k)p_{t -t_k+t_1}(x_k,x_1)\frac{\ex^{\mu  t}}{t}\,\d t .
\end{split}
\end{equation} 
Assume that $ \kappa(x)-\mu> 0 $ for at least one $ x\in\gb $. The mass of the event there is at least one jump,  $\{\exists \, t>0\colon X_t\neq X_0\} $, reads as
\begin{equation}\label{mass_boson}
M_{\gb,\mu,\beta}^B[\{ \text{ there is at least one jump } \}]=\log\left(\frac{ \det(I-\ex^{\beta(D+\mu I)}) }{\det(I-\ex^{\beta(Q+\mu I)})}\right)\, 
\end{equation}
for the Bosonic loop measure and as
\begin{equation}\label{mass_markov}
M_{\gb,\mu}[\{ \text{ there is at least one jump } \} ]=\log \left(\frac{ \det(D+\mu I) }{\det( Q+\mu I)}\right)\, 
\end{equation}
for the Markovian one. Here, $D$ is the diagonal part of the generator matrix $Q$.  The total mass of the Bosonic loop measure is given by
\begin{equation}\label{mass_boson_total}
M_{\gb,\mu,\beta}^B[\Gamma]=-\log(\det(I-\ex^{\beta(Q+\mu I)}))\, ,
\end{equation}
whereas for the Markovian loop measure we have $M_{\gb,\mu}[\Gamma]=\infty$  independent of the value of $\mu$ or the killing vector $\kappa$. For vanishing killing $\kappa\equiv 0$ and $\mu=0$ both measures have infinite mass.
\end{proposition}

\begin{remark}[\textbf{Temperature limit of the Bosonic loop measure}]
We consider the infinite temperature limit of the Bosonic loop measure with $ \mu\le 0 $ and $ \beta>0 $, that is, we replace the inverse temperature $ \beta $ by $ \beta_N=\frac{1}{N} $. 
By the definition of the Riemann integral,  it follows that 
\begin{equation}
\lim_{N\to\infty}M^B_{\Gcal,\beta_N,\mu}[G]=M_{\Gcal,\mu}[G]\, 
\end{equation}
for any $G\subset\Gamma$ measurable such that the function $t\mapsto\P_{x,x}^{\ssup{t}}(G)$ is Riemann-integrable. This allows us to interpret the Bosonic loop measure as a finite temperature version of the Markovian one.\hfill $ \diamond$ 
\end{remark}

We now fix $\gb$ to be either a finite connected subset of the integer lattice, $\gb=\L\Subset\Z^d$,   or be a finite connected subset $\L\Subset\Z^d$ times the discrete, one-dimensional torus $ \T_N=\Z/N\Z $ of length $N$, denoted by $\Lambda\times\T_N$. The torus $ \T_N $ may be represented by the lattice set $ \{0,1,\dots,N-1\} $ of size $ N$, once it is equipped with the metric $ \rho(x,y)=\inf\{\abs{x-y}_\infty\colon k\in N\Z\} $ where $ \abs{\cdot}_\infty $ denotes the maximum norm in $ \R $. In the latter case, we will refer to all objects as \textit{space-time}, e.g., $M_{\Lambda\times\T_N,\mu}$ is called the space-time Markovian loop measure with chemical potential $ \mu $. 

\medskip

\begin{definition} [\textbf{Independent space-time random walk}]\label{definition of the walk}
Let  $\Lambda\Subset\Z^d$ be finite, $\mu\le 0,\,\beta>0 $, fix a set of weights $(w(x,y))_{x,y\in\Lambda}$ and killing $(\kappa(x))_{x\in\Lambda}$ , or,  equivalently, a generator $Q$. We will refer to the space-time random walk on $ \L\times\T_N $ with generator $Q_N=Q\oplus N\beta^{-1}(\Sigma-I)$ as the \textsf{independent} case. Here $\oplus$ refers to the Kronecker plus \footnote{For $A\in\C^{n\times n} $ and $ B\in\C^{m\times m} $, define  $ A\oplus B=A\otimes \id_{\C^m}+\id_{\C^n}\otimes B$}  and $\Sigma$ is the right-shift by one on the torus, i.e.,  
$$\Sigma(\tau,\sigma)=\begin{cases} 1 & \mbox{ if } \sigma=\tau +1 \,\text{ and } \tau=0,1,\ldots, N-2,\\
1 & \mbox{ if } \tau=N-1, \sigma=0,\\
0 & \mbox{ otherwise. }
\end{cases}
$$ Equivalently, the  weights of $Q_N $ are defined as 
\begin{equation}
w_{N}(x,\tau;y,\sigma)=\begin{cases}
\beta^{-1}N\Sigma(\tau,\sigma) &\text{ if }x=y\, ,\tau\not= \sigma,\\
w({x,y}) &\text{ if }\tau=\sigma\,,x\not= y,\\
0 &\text{ otherwise.}
\end{cases}
\end{equation}
In Figure~\ref{fig1} we sketch the independent space-time random walk. In Example~\ref{counterexample} we consider the space-time random walk with symmetrised weights along the time torus. 
\begin{figure}[H]
\includegraphics[scale=0.59]{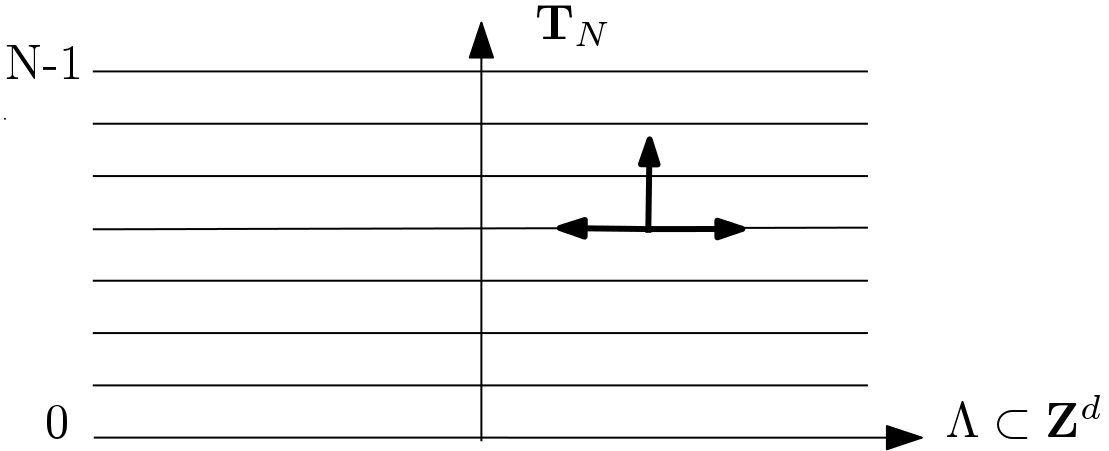}
\caption{The graphical representation of an independent space-time random walk; horizontal arrows refer to the random walk on the spatial graph $ \L $ whereas the vertical arrows indicate the walk along the ``time"-torus upwards with both walks being independent of each other.}\label{fig1}\end{figure}
\end{definition}
\begin{remark}
Note that the space-time random walk defined above winds around the torus in the "positive" direction. Equivalently, one  can use $\Sigma^{-1}$ instead of $\Sigma$ in the above definition, reversing the direction on the torus. 
In Example \ref{counterexample}  we introduce  space-time random walks with symmetric jump rates along the torus to which our results do not apply.
\end{remark}

For the following results we need to discriminate between  different classes of space-time random walks which all include the independent space-time random walk.

\medskip

\noindent \textbf{Strongly asymptotically independent space-time random walks.}\\[1ex]
Let $ G_N=Q_N+E_N $ be the generator of a space-time random walk on $ \L\times\T_N $ with  $ Q_N$ being the  generator of the independent space-time random walk (as defined above) and where the positive perturbation $ E_N $ satisfies 
\begin{equation}\label{class1}
\norm{E_N(x,y)}_1=o(\ex^{-\alpha N^{2}})  \mbox{ for some }\alpha>\beta^{-1}\max_{x\in\L} \lambda(x) \mbox{ uniformly for any }x,y\in\L, \mbox{ as } N\to \infty\, , 
\end{equation}
with $\norm{E_N(x,y)}_1 $ being the matrix norm given by the sum of the absolute values of all matrix entries of the $\T_N \times \T_N$ matrix obtained from  $E_N$ when the spatial variables $x$  and $ y $ are fixed.  A space-time random walk with generator $G_N$ satisfying \eqref{class1} is called \textit{strongly asymptotically independent space-time random walk}. 
\medskip

\noindent In the second class the decay of the perturbation is relaxed.\\[1ex]
\noindent \textbf{Weakly asymptotically independent space-time random walks.}\\[1ex]
Let $ G_N=Q_N+E_N $ be the generator of a space-time random walk on $ \L\times\T_N $ with $ Q_N$ being  the generator of the independent space-time random walk. Suppose  the positive perturbation $ E_N $ satisfies
\begin{equation}\label{class2}
\norm{E_N}_1=o\left(N^{-2N|\Lambda|}\right)\,\mbox{ as } N\to\infty,
\end{equation}
where the norm is the sum of the absolute values of the matrix entries of $E_N$ as a matrix with rows and columns indexed by $\L \times \T_N$.
We  call the space-time random walk with generator $G_N$ satisfying \eqref{class2} \textit{weakly asymptotically independent space-time random walk}.

\medskip

We will later define a more general class of space-time random walks, see Theorem~\ref{support_prop}.\\

For any path 
$$
X\colon [0,\infty)\to\L\times\T_N, t\mapsto X_t=(X^{\ssup{1}}_t,X^{\ssup{2}}_t)  
$$
we define the projection onto the spatial lattice component of the trajectory as
$$ 
\begin{aligned}
\pi_\L(X)&\colon[0,\infty)\to \L, X\mapsto\pi_\L(X) \\
& \pi_\L(X)_t:=X^{\ssup{1}}_t.
\end{aligned}
$$
We now prove the convergence of the projected space-time Markovian loop measure 
$$M_N:=[M_{\Lambda\times\T_N,\mu}]\circ\pi_\L^{-1} 
$$ to $M_{\Lambda,\mu,\beta}^B$  in the torus limit $ N\to\infty $.  The measure $M_N $ depends on the parameter $ \beta $ of the loop measure $ M_{\L\times\T_N,\mu} $ of the space-time random walk. Denote $ \bGamma $ the space of space-time loops, then $ \pi_\L\colon\bGamma\to\Gamma $, and for a measurable set $G\subset\Gamma$ we thus have
\begin{equation}\label{M_N def}
M_N[G]=\sum_{\tau=0}^{N-1}\sum_{x\in\Lambda}\int_0^\infty \frac{\ex^{\mu t}}{t}\,\P_{(x,\tau),(x,\tau)}^{\ssup{t}}(\{X=(X_s)_{s\in[0,t]}: \pi_\Lambda (X) \in G\})\d t\,.
\end{equation} 
\begin{theorem}\label{THM-main1}
Let $ \mu\le 0 $ and $ \beta>0 $. Assume that $M_N$ is induced by a strongly asymptotically independent space-time random walk with $\kappa(x)-\mu>0$ for at least one $x\in\Lambda$. For  $ k\in\N $,  let $x_1,\ldots,x_k\in\Lambda$ and $0<t_1<\cdots<t_k<\infty $ and  $A\subset (t_k,\infty)$ measurable with $\beta j\cap \partial A=\emptyset$ for all $j\in\N$. We write  $p^{\ssup{k}}:=\prod_{j=1}^{k-1}p_{t_{j+1}-t_j}(x_j,x_{j+1})$ for the jump probabilities. Then 
\begin{equation}\label{Loop_measure_convergence_theorem}
\begin{split}
& M_{N}\left[X_{t_1}=x_1,\ldots,X_{t_k}=x_k,\,\ell\in A\right]\; \underset{N\to\infty}{\longrightarrow} \beta\sum_{j=1}^\infty p^{\ssup{k}}\mathbbm{1}_A(j\beta)p_{j\beta-t_k+t_1}(x_k,x_1)\frac{\ex^{\mu j\beta}}{j\beta}\\
& = M_{\Lambda,\mu,\beta}^B\left[X_{t_1}=x_1,\ldots,X_{t_k}=x_k,\,\ell\in A\right] \, .
\end{split}
\end{equation}
Moreover, the sequence of measures is tight and we have that $M_N$ converges weakly to $M_{\Lambda,\mu,\beta}^B$ on all loops whose lengths are bounded away from zero.
\end{theorem}

\begin{remark}
The condition that $\beta j\cap \partial A=\emptyset$ for all $ j\in\N $  can be seen as an analogue to the condition in Portmanteau's theorem: Given we know $\nu_N\to\nu $ as $ N\to \infty $ weakly, one needs $\nu(\partial A)=0$ to conclude that $\nu_N(A)\to \nu(A)$ as $ N\to\infty $ for a given set $A$. 
The requirement that $A$ is bounded away from 0 stems from the fact that a loop of zero length induces a singularity in the Markovian loop measure that makes it impossible to define finite dimensional distributions.\hfill $\diamond$
\end{remark}

The following theorem is an extension of Theorem~\ref{THM-main1} by considering now  the distribution of the local time and of the occupation field. In addition we consider Poisson gases of loops defined by our Markovian and Bosonic loop measures.  The local time of a given loop $\omega $ at a point $z\in\L$ is given by
\begin{equation}
L_z(\omega)=\int_0^{\ell(\omega)}\mathbbm{1}\{X_s(\omega)=z\}\d s\, .
\end{equation}
The local time over the space-time graph $ \L\times\T_N $ is denoted by $\bL=\bL(\bomega) $ with $ \bL\colon  \L\times\T_N  \to[0,\infty) $ where $ \bL_{(x,\tau)}=\int_0^{\ell(\bomega)}\,\1\{X_s(\bomega)=(x,\tau)\}\d s,\, \bomega\in\Gamma $. Given the local time vector $ \bL $ on $ \L\times\T_N $, we obtain the local time vector $ L=(L_x)_{x\in\L} $ on $ \L $ by projection, that is, $ L=\pi_\L(\bL) $ with
$$
L_x=\sum_{\tau=0}^{N-1} \bL_{(x,\tau)},\quad x\in\L.
$$
We introduce the set-up for Poisson point measures on the set $ \Gamma $ of (rooted) loops  with $ \sigma$-algebra $ \bG $. A pure point measure $ \eta $ on $ (\Gamma,\bG) $ is a $\sigma$-finite measure of the form $ \eta=\sum_{i\in I}\delta_{\omega_i} $, where $ (\omega_i)_{i\in I} $ is an at most countable collection of loops such that $ \eta(A) <\infty $ for all $A=\{\omega\in \Gamma\colon a\le\ell(\omega)\le b\}, 0<a<b<\infty$, where
$$
\eta(A)=\#\{i\in I\colon \omega_i\in A\}, \:A\in\bG.
$$
We let $ \Omega$ be the set of pure point measures on $ (\Gamma,\bG) $, and endow $ \Omega$ with the the $ \sigma$-algebra $ \Acal $ generated by all evaluation maps $\eta\in\Omega\to\eta(A)\in\N\cup\{\infty\} $ (for details see \cite{Sznitman,LeJan}). Any $ \sigma$-finite measure on $ (\Gamma,\bG) $ defines a probability measure on $ (\Omega,\Acal) $, called the Poisson point measure, e.g., we denote by $\Poi_{\gb,\mu}, \Poi_N $ and $ \Poi_{\gb,\mu,\beta}^B $ respectively the Poisson point measures with intensity $ M_{\gb,\mu},M_N $, and $ M_{\gb,\mu,\beta}^B $, and we write $ \Er_{\gb,\mu}, \Er_{N} $ and $ \Er_{\gb,\mu,\beta}^B $ for their expectations.  Finally, for $ \eta\in\Omega, x\in\gb $, we define the occupation field of $ \eta $ at $x$ via:
\begin{equation}\label{occfield}
\begin{aligned}
\Lscr_x(\eta)&=\langle\eta,L_x\rangle\in [0,\infty],\\
&=\sum_{i\in I}L_x(\omega_i), \mbox{ if } \eta=\sum_{i\in I}\delta_{\omega_i}\in\Omega.
\end{aligned}
\end{equation}
We write $ \Lscr=(\Lscr_x)_{x\in\L} $ for the field associated to $ L $, respectively for the space-time version we write $ \bLscr $. The following theorem characterises the local time under the torus limit ($N\to\infty $) of the projected space-time loop measure.

\begin{theorem}[\textbf{Convergence of local time distributions}]\label{THM-main2}
Let $ \mu\le 0 $ and $ \beta>0 $. Assume  that the Markovian loop measure $ M_{\L\times\T_N,\mu} $ on $ \L\times\T_N $ is induced by a weakly asymptotically independent space-time random walk.
\begin{enumerate}
\item[(a)] 
Let $F:[0,\infty)^\Lambda\to\R$ be continuous and bounded such that $F(0)=0$ and the (right) derivative $\partial_{x} F$ at $0$  exists for all $ x\in\L $. Then
\begin{equation}\label{Con}
M_N[F(L)]=M_{\L\times\T_N,\mu}[F(\pi_\L(\bL))]\;\underset{N\to\infty}{\longrightarrow}\quad \beta \sum_{x\in\Lambda} \partial_x F(0)+M_{\Lambda,\mu,\beta}^B[F(L)]\, ,
\end{equation}
where the first equality is the definition of $M_N$, see \eqref{M_N def}.

\item[(b)] We can split the local time,
$$\bL=\bL_{\L,\T_N}+\bL_{\T_N}^\L+\bL_{\L}^{\T_N}+\bL^{\L,\T_N}\;\mbox{ and }  \;  L=L_\L+L^\L,
$$ where the subscript means a jump on this space and the superscript implies no jump, i.e.
\begin{equation}
    \begin{split}
        &\bL_{\L,\T_N}=\bL\1\{\exists\, t_1\colon X^\ssup{1}_{t_1}\neq X^\ssup{1}_{0}\text{ and }\exists\, t_2\colon X^\ssup{2}_{t_2}\neq X^\ssup{2}_{0}\}\, ,\\[1ex]
        &\bL_{\T_N}^\L=\bL\1\{\forall\, t_1\le \ell(\omega)\colon X^\ssup{1}_{t_1}= X^\ssup{1}_{0}\text{ and }\exists\, t_2\colon X^\ssup{2}_{t_2}\neq X^\ssup{2}_{0}\}\, ,\\[1ex]
       &\bL_{\L}^{\T_N} =\bL\1\{\exists\, t_1\colon X^\ssup{1}_{t_1}\neq X^\ssup{1}_{0}\text{ and }\forall\le \ell(\omega)\, t_2\colon X^\ssup{2}_{t_2}= X^\ssup{2}_{0}\}\, ,\\[1ex]
        &\bL^{\L,\T_N}=\bL\1\{\forall\, t_1\le \ell(\omega)\colon X^\ssup{1}_{t_1}= X^\ssup{1}_{0}\text{ and }\forall\, t_2\le \ell(\omega)\colon X^\ssup{2}_{t_2}= X^\ssup{2}_{0}\}\, ,
    \end{split}
\end{equation}
and similarly for $L_\L $ and $ L^\L$, respectively.
If the Markovian loop measure is induced by the independent space-time random walk, then as $N\to\infty$,   
\begin{equation}\label{convergence split}
\begin{aligned}
& M_{\Lambda\times \T_N}[F(\pi_\L(\bL_{\L,\T_N}))]\longrightarrow M^B_{\Lambda,\mu,\beta}[F(L_\Lambda)]\,,\\[1ex]
& M_{\Lambda\times \T_N}[F(\pi_\L(\bL_{\T_N}^\L))]\longrightarrow  M^B_{\Lambda,\mu,\beta}[F(L^\Lambda)]\, ,\\[1ex]
& M_{\Lambda\times \T_N}[F(\pi_\L(\bL_{\L}^{\T_N}))]\longrightarrow   0  \, ,\\[1ex]
&  M_{\Lambda\times \T_N}[F(\pi_\L(\bL^{\L,\T_N}))]\longrightarrow   \beta\sum_{x\in\Lambda}\partial_xF(0)\, ,
\end{aligned}
\end{equation}
with the same assumptions on $F$ as in $(a)$.
\item[(c)]  Denote $\Er_{N}$ the expectation with respect to the Poisson Point Process associated with $M_N$. Then, for every bounded and continuous function $F\colon[0,\infty)^\Lambda\to\R$ we have that
\begin{equation}
\Er_{N}[F(\Lscr)]\longrightarrow \Er_{\L,\mu,\beta}^B[F(\Lscr+\beta)]\,\mbox{ as } N\to\infty\, ,
\end{equation}
with the results from \eqref{convergence split} carrying over to the occupation field.
Here, the constant occupation field $\beta$ represents the mean time it takes for a space-time random walk to wind around the torus once.
\end{enumerate}
\end{theorem}

We return to the setting of a general set $\gb$ endowed with the graph structure inherited by the weights. The splitting of the local times  can be analysed in this setting as well.  We call loops with no jumps (loops are not leaving  their origin) on $ \gb $ \textit{point loops} and denote their corresponding local time vector by $ L^{\gb} $ and their  occupation field by  $ \Lscr^{\gb} $ whereas $ L_{\gb} $ and $ \Lscr_{\gb} $ are the corresponding functionals for \textit{genuine loops}, i.e., loops which are not point loops,  (the role of the sub- and superscripts corresponds to the one defined above).

\begin{proposition}[\textbf{Local times}]\label{diagonalization_theorem}
Let $ \beta> 0, \mu\le 0 $, and   denote  $(d(x))_{x\in\gb}$  the diagonal elements of $Q+\mu I$.  The occupation field can be written as a sum of two parts, $\Lscr=\Lscr^{\gb}+\Lscr_{\gb}$, the point and the genuine loops, with the following properties:
\begin{enumerate}
\item[(a)]
The splitting is independent, i.e.,  $\Lscr^{\gb}\bot\; \Lscr_{\gb}$.

\medskip

\item[(b)]
Under the Poisson measure $\Pr_{\Lambda,\mu}$ the occupation field of point loops $\Lscr^{\gb}$  factorises into Gamma distributed random variables with parameters $(1, d(x))_{x\in\gb}$. Under the Bosonic Poisson measure $\Pr_{\gb,\mu,\beta}$, the occupation field of point loops $\Lscr^{\gb}$ is distributed like a product of random variables with Laplace transform 
$$
\Big(\E[\ex^{-v(x)\Lscr^{\gb}_x}]\Big)_{x\in\gb}=\Big(1-\frac{\ex^{\beta(d(x)+\mu)}\left(\ex^{-\beta v(x)}-1\right)}{1-\ex^{\beta(d(x)+\mu)}}\Big)_{x\in\gb}\,\quad \mbox{ with  } \; v(x)\in\R_+, x\in\gb\,.
$$

\end{enumerate}
\end{proposition}

The following example shows that symmetric random walk on the torus does not generate Bosonic loops by spatial projection.


\begin{figure}[H]
\includegraphics[scale=0.58]{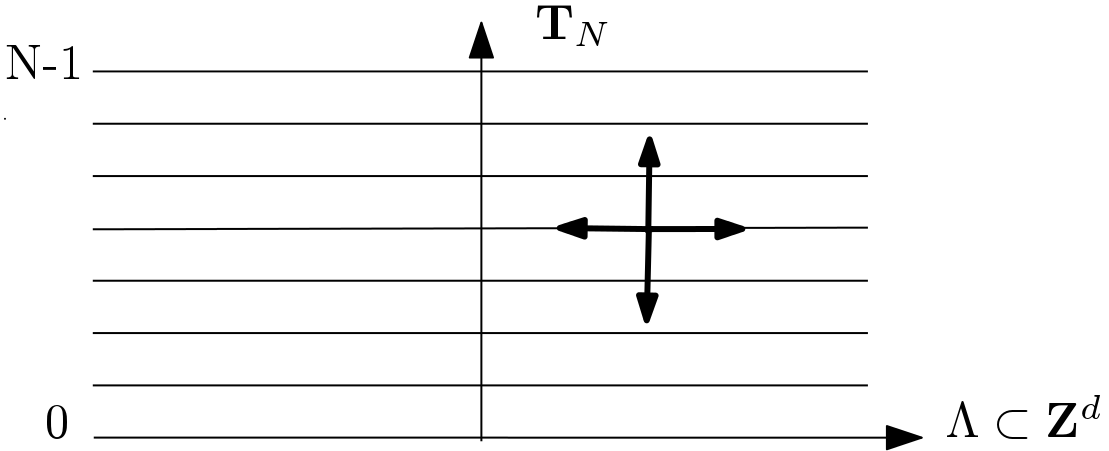}
\caption{The graphical representation of a space-time random walk for which the above convergence results are no longer valid}\label{counterexample_picture}
\end{figure}
\begin{example}\label{counterexample}
We consider the``symmetrisation" of the torus walk of the space-time random walk described above (see Figure~\ref{counterexample_picture}). We go up or down on the tours independently with rate $\beta^{-1}N$ and obtain the space-time rates
\begin{equation}
\tilde{w}_{\gb}(x,\tau;y,\sigma)=\begin{cases}
\beta^{-1}N(\Sigma(\tau,\sigma)+\Sigma(\sigma,\tau))&\text{ if }x=y,\\
w({x,y})&\text{ if }\sigma=\tau,\\
0 &\text{ otherwise}
\end{cases}
\end{equation}
on $ \gb=\L\times\T_N $. This space-time random walk has a symmetric generator but  Theorems \ref{THM-main1} and \ref{THM-main2} fail to hold. The proof is given in the Section \ref{section proofs}. \hfill $ \clubsuit$ 
\end{example}

\begin{figure}[H]
\includegraphics[scale=0.6]{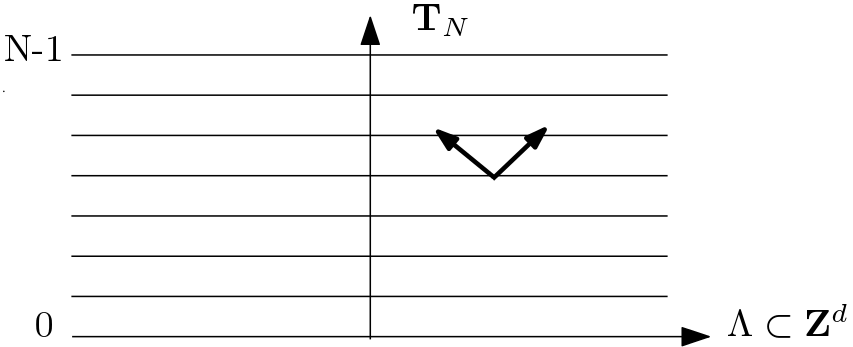}
\caption{A  space-time random walk where torus and lattice jumps are correlated}\label{mixingpic}\end{figure}

Finally, the third class is related to Bosonic loops and is not defined via perturbations of the independent space-time random walk.\\[1ex]
\noindent \textbf{Bosonic Random Walk.}\\[1ex]
Let $ (G_N)_{N\in\N}$ be a family of generators on $\Lambda\times\T_N$ with weights $(w_N (v,u))_{v,u\in\Lambda\times\T_N}$ and $\kappa_N(u)-\mu$ uniformly (in $N$) bounded away from 0 for at least one $ u\in\L\times\T_N $. Let $z=(z_1,\ldots,z_n)\in (\Lambda\times\T_N)^n$ be a path. The coordinate projections for $ z_i=(z_i^{\L},z_i^{\T_N}) $ are $ z_i^{\L}\in\L $ and  $z_i^{\T_N}\in\T_N$ respectively. For any such path $z$, denote $\wind(z)=\#\{i>1\colon z_i^{\T_N}=z_1^{\T_N}\}$ the torus winding number for the given path. If $z_1=(x,\tau)$ and $z_n=(y,\sigma)$ write $z\colon (x,\tau)\to(y,\sigma)$.\\
For such a path $z$ denote 
\begin{equation}\label{step probabolities}
P(z)=\prod_{i=1}^{n-1}p(z_i,z_{i+1})\, ,
\end{equation}
the product of the jump probabilities, and denote $ W(z) $ the random  waiting time which is distributed as the sum $ \sum_{i=1}^n X_i $ of independent exponentially distributed random variables $ X_i\sim\Ex((-G_N(z_i,z_i))) $.

Assume that for all $x,y\in\Lambda$ there exists  $d_{x,y}\in (0,1]$ such that
\begin{equation}\label{waitprob}
\lim_{N\to\infty}\frac{1}{N}\sum_{\tau=0}^{N-1}\sum_{\heap{z\colon (x,\tau)\to(y,\tau)}{\wind(z)=1}}P(z)=d_{x,y}\, ,
\end{equation}
and
\begin{equation}\label{waitingtime}
\lim_{N\to\infty}\frac{1}{N}\sum_{\tau=0}^{N-1}\sum_{\heap{z\colon (x,\tau)\to(y,\tau),}{\wind(z)=1}}P(z)\big[\big(\beta-W(z)\big)^2\big]=0\, .
\end{equation}

%
Finally, assume that $\tau\in\T_N$,  we have $-G_N(u,u)=N\beta^{-1}+o(N)$ for any  $ u=(x,\tau)\in\L\times\T_N$.\\
If all the above conditions are satisfied we say that the space-time random walk with generator $G_N $ is a \textit{Bosonic random walk}. Note that \eqref{waitingtime} ensures that in the torus limit the turnaround time for the torus is precisely $ \beta $.

We will show that the length of any loop of the Bosonic random walk converges to a multiple of $\beta$ in a suitable way, that is, the spatial projections of their loop measures have restricted support, namely they are supported on the Bosonic time horizons.

\begin{theorem}{\textbf{Convergence of the support. }}\label{support_prop}
For any $ N\in\N $ let $ G_N $, be the generator of a  Bosonic random walk, and denote $M_N $ the spatial projection of the loop measure onto $ \L $. Let  $A\subset (0,\infty)$ be measurable, bounded away from zero and bounded from above. Then 
\begin{equation}
\lim_{N\to\infty}M_N\left[\ell\in A\right]
>0 \,\,\text{   if   }\,\, \exists k\in\N\colon k\beta\in A^\circ\, ,
\end{equation}
and
\begin{equation}
\lim_{N\to\infty}M_N\left[\ell\in A\right]
=0\,\, \text{   if   }\,\, \forall k\in\N\colon k\beta\in \left(\R\setminus A\right)^\circ\, .
\end{equation}
Furthermore, all weakly and strongly asymptotically independent space-time random walks with spatial killing such that $ \kappa(x)-\mu >0 $ for at least for one $x\in\L $, are Bosonic random walks.
\end{theorem}

The following example of a Bosonic walk shows an interesting phenomenon when we mix the spatial component with the torus component of the space-time random walk, see Figure~\ref{mixingpic}.
\begin{example}\label{ergodic example}
Let $G_N$ be the space-time random walk on $\gb=\L_M\times\T_N $ with $\L_M=\{-M,\ldots,M\}^d\Subset\Z^d$ for some $M>0$ equipped with periodic boundary conditions, i.e.,  the  rates are 
$$
w_N^{\texttt{per}}(x,\tau,y,\sigma)=\begin{cases}
N\beta^{-1}/2d&\text{ if }\sigma=\tau+1 \mod(N),\, x\sim y\, ,\\
0&\text{ otherwise.}
\end{cases}
$$
Here, by $x\sim y$ we mean that $x_i-y_i=0$ for all but one $i\in\{1,\ldots, d\}$, and for this specific one, say  $i_0$, we have $x_{i_0}-y_{i_0}=1\mod (2M+1)$. For $ k\in \N $,  let $  x_i\in\L $ for $ i=1,\ldots, k  $,  be given and let $ 0<t_1<t_2<\cdots< t_k $. 
Let $ A\subset\R_+ $ be measurable with $ A\subset(t_k,\infty) $ and $ \beta n\cap \partial A=\emptyset $ for all $ n\in\N $.  Denote $\bpi$ the product $\prod_{i=1}^k \pi(x_i)$ where $\pi$ is the stationary distribution associated to the jump-chain on $\Lambda_M$ with weight $w(x,y)=\1\{x\sim y\}$. Then, the finite dimensional distributions $ $ from Theorem \ref{THM-main1} converge as follows.
\begin{equation}\label{example-equ-erg}
M_N\big[X_{t_1}=x_1,\ldots, X_{t_k}=x_k,\ell\in A\big]\underset{N\to\infty}{\longrightarrow} \bpi\sum_{j=1}^\infty\mathbbm{1}_A(j\beta)\frac{\ex^{\mu j\beta}}{j}\, ,
\end{equation}
which can be interpreted as Bosonic loop measure in stationary distribution. This notion is not rigorous, as the dependence on $t_1,\ldots,t_k$ is lost due to the mixing. For a sketch of the  proof, see Section~\ref{Proof1}. \hfill $ \clubsuit$ 
\end{example}
\subsection{Isomorphism Theorems and space-time correlation functions}\label{Isomorphism Theorems and space-time correlation functions}
Before we study our new isomorphism theorems  and moment formulae in Section~\ref{S-ISO},  we give some succinct  background in Section~\ref{backgroundISO}. In Section~\ref{S-space-time} we return to our previous space-time setting and demonstrate that tori limit of moments lead to probabilistic representations of quantum correlation functions.

\subsubsection{Background Isomorphism theorems and complex Gaussian measures} \label{backgroundISO}
The Markovian loop soup owes its conception to the physics community, where it arises via a functional integral description of a lattice model. In \cite{Sym69}, Symanzik provided a heuristic description of $ \phi^4$-quantum field theory in terms of a gas of interacting Brownian loops.   On considering lattice field theories in place of Symanzik's continuum model, Brydges, Fr\"ohlich and Spencer \cite{BFS82} were able to make rigorous the connection between the two models. A version of this connection can easily be seen in the random walk representation of the two-point correlation function for Gaussian fields over the integer lattice given by the Green function of the simple random walk on $ \Z^d$. Inspired by the work of Symanzik and Brydges \textit{et al.}, Dynkin \cite{Dyn84} provided an extension to continuous time processes. Dynkin's isomorphism theorem relates the local times distribution to the distribution of squares of the Gaussian field (\cite{Sznitman}).  
Loop measures have been derived in \cite{B92} via weights for discrete time random walks, and if one wishes to consider directed walks (edges) in the underlying graph, Brydges shows in \cite{B92} that we need to consider complex-valued random fields instead of only real-valued ones. 
We consider standard complex integration in $|\gb|$ variables denoted by $\prod_{x\in\gb}\d \phi_x=\prod_{x\in\gb}\d\Re(\phi_x)\d\Im(\phi_x) $, where $ \Re(\phi_x) $ is the real and $ \Im(\phi_x) $ is the imaginary part of $ \phi_x$.
%
%
%

For the convenience of the reader we cite the following theorem in \cite{BIS09}. 
\begin{theorem}[\textbf{Complex Gaussian, \cite{BIS09}}]\label{Gaussian_det_theorem}
Let $C\in \C^{\gb\times\gb}$ be a matrix with positive  definite Hermitian part, i.e., $ \langle \phi,(C+C^*)\bar{\phi}\rangle=\sum_{x,y} \phi_x(C({x,y})+\overline{C}({y,x}))\bar{\phi}_y > 0 $ for all nonzero $ \phi\in\C^\gb $, and inverse $ A=C^{-1} $. Then we define the complex Gaussian measure $ \mu_A $ on $ \C^{\gb} $ with covariance $C$, namely
$$
\d\mu_A(\phi)=\frac{1}{Z_A}\ex^{-\langle \phi,A\bar{\phi}\rangle}\,\prod_{x\in\gb}\d\phi_x,
$$
and normalisation
\begin{equation}\label{detid}
Z_A=\int\, \ex^{-\langle \phi,A\overline{\phi}\rangle}\,  \prod_{x\in\gb}\d\phi_x  =\frac{(\pi )^{|\gb|}}{\det A}\, .
\end{equation}
\end{theorem}
We write $ \E_A $ for the expectation with respect to the measure $ \mu_A $ in the following.
\begin{theorem}[\textbf{Dynkin's Isomorphism, \cite{Dyn84}, \cite{B92}}]\label{dynkin}
For any bounded measurable $ F\colon\R^\gb\to\R $
$$
\E_{x,y}^{0,1}\otimes \E_{A}\Big[F\big(L+|\phi|^2\big)\Big]=\E_{A}\Big[\phi_x\phi_y F\big(|\phi|^2\big)\Big]\, ,
$$
where $ \E_{x,y}^{0,1} $ denotes integration with  respect to $  \P_{xy}^{0,1}:=\int_0^\infty \,\P^{\ssup{t}}_{xy}\,\d t  $ and $|\phi|^2=(\overline{\phi}_x\phi_x)_{x\in\gb}$.
\end{theorem}
As observed by Le Jan \cite{LeJan}, there is a relation between the field of occupation times for the Poisson loop process with Poisson measure $ \Pr^{}_{\gb,\mu\equiv 0} $ and intensity measure $ M_{\gb,\mu\equiv 0} $ and the Gaussian measure.
\begin{theorem}[\textbf{Le Jan's Isomorphism, \cite{LeJan}}]  Let the generator $Q$ on the finite graph $ \gb $ be symmetric and  $ A=-Q $,  and denote $ \mu_A $ the unique Gaussian measure whose covariance is the corresponding Green function and whose Dirichlet form is given as $ \langle\phi, A\phi\rangle $. Then 
$$
(\Lscr_x)_{x\in\gb} \;\mbox{ under }\; \Pr^{}_{\gb,\mu\equiv 0}\;,\mbox{ has the same law as }\; (\phi_x^2)_{x\in\gb} \;\mbox{ under }\; \mu_A.
$$
\end{theorem}


\subsubsection{Isomorphism theorems}\label{S-ISO}
We shall study not only the square of the field, but also the full field in its generality. This is an important extension, but we need another input, namely to consider non-symmetric generator matrices (\cite{LJ07}) as the generators of our  space-time random walks are non-symmetric.
We are interested in deriving corresponding isomorphism theorems, first for our Markov loop measures having chemical potential $ \mu\le 0 $, and then secondly for the Bosonic loop measures with $ \beta>0, \mu\le 0$. The distribution of the occupation field for Bosonic loop measures is not given by a Gaussian measure - in fact they are related to Permanental processes (cf. \cite{FR14}). \\

The next theorem is a natural extension of Le Jan's loop soup version of  Dynkin's Isomorphism theorem (see \cite{LJ07},\cite{B92}).
\begin{theorem}[\textbf{Isomorphism theorems}]\label{Isomorphism occupation field}
\begin{enumerate}
Let $A=-(Q+\mu\1)$, where $Q$ is the generator  on the finite graph $ \gb$ such that $ \kappa(x)- \mu> 0 $ for at least one $x\in\gb$ for  $ \mu\le 0$.  Denote $ \mu_A $  the complex Gaussian measure defined in Section~\ref{backgroundISO}. We write $ \Er_{\gb,\mu} $ for the expectation with  respect  to  the Poisson process with intensity measure being the Markovian loop measure $ M_{\gb,\mu} $. Then the  following holds.
\item[(a)]
The measure  $\mu_{A}$ evaluated on the squares $ |\phi|^2=(\overline{\phi}_x\phi_x)_{x\in\gb} $ of the random field is a probability measure. Furthermore, for every bounded and continuous  function $F\colon \R_+^\gb\to \C$,
\begin{equation}\label{Isoa}
\E_A[F(|\phi|^2)]=\Er_{\gb,\mu}[F(\Lscr)]\, .
\end{equation}

\medskip

\item[(b)]
Denote $\widetilde{\Er}_{\gb,\mu}$ the expectation with respect to the Poisson process with intensity measure being the Markovian loop measure for the random walk with  the symmetrised generator $\widetilde{Q}:=\frac{1}{2}(Q+Q^T)$. Then, for every $F\colon\C^\gb\to\C$ bounded and continuous, 
\begin{equation}
\E_A[F(\phi)]=\widetilde{\Er}_{\gb,\mu}\Big[\frac{1}{Z_S\int\ex^{-\langle \sqrt{\Lcal}\theta,\widetilde{A}\sqrt{\Lcal}\bar{\theta}\rangle}\d \Scal_\gb(\theta)}\int F(\theta\sqrt{\Lcal})\ex^{-\langle \sqrt{\Lcal}\theta,A\sqrt{\Lcal}\bar{\theta}\rangle}\d\Scal_\gb(\theta)\Big]\, ,
\end{equation}
where $Z_S=\frac{\det(\widetilde{A})}{\det(A)}$ with  $\widetilde{A}=-(\widetilde{Q}+\mu\1) $ being the Hermitian part of $A$. The measure $\d\Scal_\gb$ is defined as the product of the uniform measure on $\{z\in\C\colon |z|=1\}$ over all $x\in\gb$. 

\medskip

\item[(c)]
For every  $F\colon\C^\gb\to\C$ bounded and continuous, 
\begin{equation}
\E_A[F(\phi)]=\Er_{\gb,\mu}\Big[\frac{1}{\int \ex^{-\langle\theta\sqrt{\Lcal},A\sqrt{\Lcal}\bar{\theta}\rangle}\d \Scal_\gb(\theta)}\int F(\theta\sqrt{\Lcal})\ex^{-\langle\theta\sqrt{\Lcal},A\sqrt{\Lcal}\bar{\theta}\rangle}\d \Scal_\gb(\theta)\Big]\, .
\end{equation}
\end{enumerate}
\end{theorem}

\begin{remark}
The first identity, which  appeared in \cite{LJ07},  is more surprising than it may seem. Since $A$ may not be a symmetric matrix, it is a priori not clear that the complex-valued distribution of $\phi_x\overline{\phi_x}$ does in fact only takes positive values. This has been remarked by \cite{LJ07}, however in the context of an isomorphism theorem of the local time vector. In \cite{LeJan} a similar identity to part (a)  is mentioned  in the context of symmetric processes. The second statement is  motivated by recent studies  in \cite{Cam15} on symmetric generators. The last identity is new, to the best of our knowledge.
\hfill $\diamond $
\end{remark}

We now state and prove a generalised version of Symanzik's moment representation formula. It expresses the moments of certain interacting fields in terms of a Poisson gas of loops and a collection of paths interacting with the loops.  In Section~\ref{S-space-time}  we  apply that representation  to our space-time random walk stetting.
In what follows $ \mu_A $ denotes  the Gaussian measure  with  $ A=-(Q+\mu I) $  defined in Section~\ref{backgroundISO}. We define a class of perturbed Gaussian measures. Let $J\colon [0,\infty)^\gb\to\R$ be bounded and continuous function of $ |\phi|^2=(\overline{\phi}_x\phi_x)_{x\in\gb} $ such that
\begin{equation}\label{Jclass}
\int\,\ex^{-\langle\phi,A\overline{\phi}\rangle}J(|\phi|^2)\,\prod_{x\in\gb}\d\phi_x \in (0,\infty)
\end{equation} holds. The perturbed Gaussian measure $ \mu_{A,J} $ is defined as
\begin{equation}
\frac{\d\mu_{A,J}}{\d\mu_{\beta A}}(\phi)=\frac{Z_{\beta A}}{Z_{A,J}}J(\beta|\phi|^2)\, \mbox{ with } Z_{A,J}=\int\,\ex^{-\langle\phi,\beta A,\overline{\phi}\rangle}J(\beta |\phi|^2)\prod_{x\in\gb}\d\phi_x . 
\end{equation}
Expectation with respect to $ \mu_{A,J} $ is denoted  $\E_{A,J}$. Furthermore, for any $ \beta>0 $ and $ \mu\le 0 $, define 
\begin{equation}\label{bridgemeasure}
\P_{xy}^{\mu,\beta}[G]:=\int_0^\infty \ex^{t\beta\mu}\,\P^{\ssup{t\beta}}_{xy}[G]\,\d t, \quad \mbox{ for }G\subset D_\gb \mbox{ measurable}\,.
\end{equation}

The following is a generalisation of Symanzik's formula (see \cite{Sznitman} for the standard version for real-valued Gaussian Free Fields).

\begin{proposition}[\textbf{Moments of perturbed Gaussian}]\label{Symanzik_formula_generalized}
Let $J$ be bounded and continuous satisfying \eqref{Jclass}, and $ \beta>0 $ and $ \mu\le 0 $, and assume that $ \kappa(x)-\mu>0 $ for at least one $ x\in\gb$.  Then
\begin{equation}\label{moment}
\E_{A,J}\big[\bar{\phi}_x\phi_y\big]=\frac{\E_{xy}^{\mu,\beta}\otimes\Er_{\gb,\mu}\left[J(\Lscr+L) \right]}{\Er_{\gb,\mu}\left[J(\Lscr)\right]}\, .
\end{equation}
Here  $L\sim\P_{xy}^{\mu,\beta}$ and $\Lscr\sim\Pr_{\gb,\mu}$. Higher even moments can be obtained analogously, see \cite{Sznitman}.
\end{proposition}
\begin{remark}
This is a generalisation of Symanzik's formula in \cite{Sznitman} where $J$ is chosen to be of the form
\begin{equation}
J(u)=\prod_{x\in\gb}\int_0^\infty\ex^{-u_xw}\nu(\d w)\, ,\qquad \mbox{ with } \nu\in\Mcal_1([0,\infty))\,, u=(u_x)_{x\in\gb}\in [0,\infty)^{\gb}.
\end{equation}
Of interest are certainly $ \phi^4$-perturbations which we are going to study in the future. 
\end{remark}

\subsubsection{Moments for space-time random walks}\label{S-space-time}
We return to the space-time setting with  $\gb=\Lambda\times\T_N$ and the non-symmetric generator $G_N$ of some weakly asymptotically independent space-time random walk. Let $ A_N=-(G_N+\mu I) $ and $ A=-(Q+\mu I)$.  We aim  to apply formula \eqref{moment} with  the non-symmetric matrix $A_N$ and to perform the torus-limit $ N\to\infty $ for spatial projections. It will turn out that the right hand side of \eqref{moment} for spatial projections converges to the one-particle reduced density matrix (quantum correlation function) represented in Theorem~\ref{physical quantities representation}).  We  define spatial  projections of perturbed complex Gaussian measures $ \bmu_{A_N,\Jsf} \in\Mcal_1(\C^{\L\times\T_N})$ and interaction functionals suitable for Boson particle systems (see Section~\ref{S:Bosons}). 
Suppose that $ \mathsf{V}\colon[0,\infty)^{\L\times\T_N}\to\R $ is  continuous such that $ \mathsf{J}:=\ex^{-\mathsf{V}} $ is bounded, continuous and satisfies condition \eqref{Jclass} to ensure that the measure $ \bmu_{A_N,\Jsf} $ is well-defined, and in addition we assume that the spatial projections, denoted $ V $ and $ J $,  satisfy equivalent conditions. That is, let 
$$ 
\psf_\L\colon\C^{\L\times\T_N}\to\C^{\L}, \bphi\mapsto\psf_\L(\bphi)=\big(\sum_{\tau=0}^{N-1}\bphi_{x,\tau}\big)_{x\in\L} 
$$ be the spatial projection, then $ V=\Vsf\circ\psf_\L^{-1} $ and $ J=\Jsf\circ\psf_\L^{-1}=\ex^{-\bV} $ ensure that $ \mu_{A,J}:=\bmu_{A_N,\Jsf}\circ\psf_\L^{-1} $ is well-defined as well. The choice of our interaction functional $ V $ is due to our occupation field distribution in Proposition~\ref{Symanzik_formula_generalized} and comprises many mean-field models for systems of Bosons, see Section~\ref{S:Bosons} and \eqref{interactionfield}. In Theorem~\ref{physical quantities representation} we demonstrate that on the Poisson process level one can incorporate typical  particle interaction functionals studied in physics (e.g. \eqref{interaction}). Having set the stage for  space-time versions of the right hand side of \eqref{moment} in Proposition~\ref{Symanzik_formula_generalized}, we only need to address the spatial projection of the path measure (see \eqref{bridgemeasure}),
\begin{equation}\label{pathmeasureST}
\P_{(x,\tau),(y,\sigma)}^{\mu,\beta}=\int_0^\infty\,\ex^{t\beta\mu}\,\P^{\ssup{\beta t}}_{(x,\tau),(y,\sigma)}\,\d t\,\quad x,y\in\L,\tau,\sigma\in\T_N,
\end{equation}
which is given for any event $ G $ of the spatial path $ x\longrightarrow y $ as
\begin{equation}\label{projpathmeas}
\P_{x,y,N}^{\mu,\beta}  (G)=\sum_{\tau=0}^{N-1}\P_{(x,\tau),(x,\tau)}^{\mu,\beta}(X^{\ssup{1}}\in G),
\end{equation}
and $ \E_{x,y,N}^{\mu,\beta} $ denotes expectation with respect to $ \P_{x,y,N}^{\mu,\beta} $.
In what follows only the spatial path $ x\longrightarrow y $ appears in the spatial projection which amounts to summing over all $ \tau $ and $ \sigma$. However, as we aim to obtain the quantum correlation function in the torus limit, we consider only the case that $ \tau=\sigma$. This is due to the standard loop representation of quantum correlation functions (see \cite{BR97}), where the open spatial path $x\longrightarrow y  $ originates from breaking any possible loop. Therefore, the condition  $ \tau=\sigma$  ensures that the spatial path $ x\longrightarrow y $ comes from a cycle in the space-time setting.  We aim to study the other case in the future.  To summarise, we shall take the torus limit of the spatial projections of 
\begin{equation}\label{momentST}
\E_{A_N,\Jsf}\big[\overline{\bphi}_{(x,\tau)}\bphi_{(y,\tau)}\big]=\frac{1}{\Er_{\L\times\T_N,\mu}[\Jsf(\bLscr)]}\E_{(x,\tau),(y,\tau) }^{\mu,\beta}\otimes\Er_{\L\times\T_N,\mu}\big[\Jsf(\Lb+\bLscr)\big]\,\quad, x,y\in\L, \tau\in\T_N.
\end{equation}
If $ \Vsf$ is linear, i.e., there exists $ \bv\in[0,\infty)^{\L\times\T_N} $ such that $ \Vsf(\cdot)=\langle\bv,\cdot\rangle $, the measure $ \bmu_{A_N,\Jsf} $ is still a Gaussian measure and the left hand side in \eqref{momentST} is the second Gaussian moment which  equals the Green function $ \bG_N(x,\tau;y,\tau) $. 
In what follows we consider $ \L\Subset\Z^d $ and the simple random walk on $ \Z^d $, i.e., the generator $Q$ is the discrete Laplacian in $ \L $ with Dirichlet boundary conditions. The following results are also true in more general cases but our choice is adapted to our results in Section~\ref{S:Bosons}.    
The crucial step is that the projected path measure \eqref{projpathmeas} converges to the Bosonic path measure,
\begin{equation}\label{bosonbridgemeasure}
\P_{x,y,B}^{\mu,\beta}(G)=\sum_{j= 1}^\infty \ex^{\beta \mu j}\P_{x}[G\cap \{X_{\beta j}=y\}]\,, \; G\subset D_\L \mbox{ measurable}.
\end{equation}

In the following, let  $v\colon[0,\infty)^\L\to\R $ be a functional such that 
\begin{equation}\label{interactionfield}
V\colon\Omega\to\R,\; \eta\mapsto V(\eta):=v(\Lscr(\eta)),
\end{equation}
with $ \Lscr_x(\eta)=\sum_{i\in I}L(\omega^{\ssup{i}}) $ for $ \eta=\sum_{i\in I}\delta_{\omega^{\ssup{i}}} $.

\medskip

\begin{theorem}[\textbf{Space-time loop-representation of the correlation function}]\label{Space-time loop-representation of the correlation function}
Suppose $ \beta>0 $ and $ \mu\le 0 $. Let $ A_N=-(G_N+\mu I) $ with $G_N=Q_N+E_N$ being  the generator of the weakly asymptotically independent space-time random walk. 
\begin{enumerate}
\item[(a)] 
\begin{enumerate}
\item[(i)] For any continuous bounded $ F\colon [0,\infty)^\L $,
$$
\E_{x,y,N}^{\mu,\beta}\big[F(L)\big]\longrightarrow  \E_{x,y,B}^{\mu,\beta}\big[F(L)\big]  \;\mbox{ as } N\to\infty.
$$
\item[(ii)]  Suppose that $ \Jsf=\ex^{-\Vsf} $ is bounded and continuous and satisfies \eqref{Jclass}, and let  $J=\Jsf\circ\psf_\L^{-1}$. Let  $V_\beta(\cdot):=V(\cdot+\beta) $ with $V$ being defined above. Then, for $ x\not= y\in \L $, 
\begin{equation}\label{RDMresult1}
\begin{aligned}
 \frac{ \E_{x,y,N}^{\mu,\beta}\otimes\Er_{N}\big[J(L+\Lscr)\big] }{\Er_N[J(\Lscr)]}\;\underset{N\to\infty}{\longrightarrow}\; & \frac{\E^{\mu,\beta}_{x,y,B}\otimes\Er^B_{\L,\mu,\beta}\big[J(L+\Lscr+\beta)\big]}{\Er^B_{\L,\mu,\beta}[J(\Lscr+\beta)]}\\[1ex]& =:\rho_{\L,V_\beta}^{\ssup{1}}(x,y)\,, 
\end{aligned}
\end{equation}
where $ \rho_{\L,V_\beta}^{\ssup{1}} $ will be later defined in Section~\ref{S:Bosons} as the one-particle reduced density matrix (correlation function)  for a system of Bosons in  $ \L $ with interaction $ V_\beta$  in thermodynamic equilibrium at inverse temperature $ \beta $ and chemical potential $ \mu $ (see Theorem~\eqref{physical quantities representation}).

\end{enumerate} 
\medskip

\item[(b)]
Suppose $ \bv\in[0,\infty)^{\L\times\T_N} $ and that   $ \Vsf(\cdot)=\langle\bv,\cdot\rangle $. Then  $ V(\eta)=\langle v,\Lscr\rangle $ with $ v=\bv\circ\psf_\L^{-1}\in [0,\infty)^\Lambda$. The torus limit of the sum of space-time  Green functions  is
\begin{equation}
\sum_{\tau=0}^{N-1}\bG(x,\tau;y,\tau)= \sum_{\tau=0}^{N-1}\E_{A_N,\Jsf}\big[\overline{\bphi}_{(x,\tau)}\bphi_{(y,\tau)}\big] \;  \underset{N\to\infty}{\longrightarrow}\;\rho_{\Lambda,V}^{\ssup{1}}(x,y),\quad x\neq y\in \L. 
\end{equation}
\end{enumerate}
\end{theorem}

\begin{remark}
\begin{enumerate}
\item[(i)]  The constant background field $ \beta $ in \eqref{RDMresult1} is due to the torus point loops, and thus they are an artefact of our method. We can obtain the limit without the constant field by restricting to genuine torus loops. 

\medskip

\item[(ii)]  One can consider true (or genuine) path interaction where each point of one path interacts with any point of the other path (double time integration) instead of the particle interaction defined in \eqref{interaction}, namely, 
\begin{equation}\label{interactionP}
\begin{aligned}
V(\omega^{\ssup{1}},\ldots,\omega^{\ssup{N}})=\frac{1}{2}\sum_{1\le i,j\le N}\int_0^{\ell(\omega^{\ssup{i}})}\int_0^{\ell(\omega^{\ssup{j}})}\,v(|\omega^{\ssup{i}}(t)-\omega^{\ssup{j}}(s)|)\,\d t\d s,
\end{aligned}
\end{equation}
where $ \ell(\omega^{\ssup{i}}) $ is the length of loop $ \omega^{\ssup{i}} $. If we restrict our loop measure to loops whose lengths are  bounded away from zero and if we consider this path interaction, we  obtain  limits like in \eqref{RDMresult1} but without the constant background field $ \beta $. On the other hand, genuine path interactions (e.g.,  Polaron type interactions in \cite{DV83,ABK04,ABK05}) are not derived from quantum systems via the Feynman-Kac formula  (see Section~\ref{S:Bosons}). 

\medskip

\item[(iii)] Item (b) in Theorem~\ref{Space-time loop-representation of the correlation function} demonstrates that quantum correlation functions can be obtained as limits of space-time Gaussian moments as long as the interaction functional is linear.   
More general interactions would require computing non-Gaussian moments followed by taking the torus limit. If one  puts $ \bv\equiv 0 $, then one can consider the so-called thermodynamic limit $ \L\uparrow\Z^d $ to obtain the one-particle reduced density matrix (correlation function) in \eqref{RDM}, see Section~\ref{sec-BEC}.  \hfill$\diamond$
\end{enumerate}
\end{remark}

\subsection{The infinite volume Limit}\label{S-limit}
We now extend our study  to countable graphs, and for ease of notation  we consider the integer lattice $\Z^d $ and the simple random walk  for ease of notation. We are interested which of our previous results carry over taking limits $ \L\uparrow\Z^d $.  An application driven motivation   is to  study the  Bosonic loop measure in the so-called thermodynamic limit $ \L\uparrow\Z^d$ , see Section~\ref{sec-bosons}.  We list some aspects  of taking the limit to $ \Z^d $ and leave a thorough study for future work.

Many of our results involve Laplace transforms and thus determinants of finite-dimensional matrices, which in the thermodynamic limit  pose technical challenges as e.g. the generator $Q$ of the simple random walk on $\Z^d$ fails to be trace-class.  On the other hand, we can expect that some of the results continue to hold. We summarise which of our previous  results  continue to hold in thermodynamic limit.

\begin{theorem}\label{infvollimprop}
The following results  hold in the thermodynamic limit $\Lambda\uparrow\Z^d$.
\begin{enumerate}
\item[(a)] Proposition~\ref{properties}: The formulae for the finite dimensional distributions carry over  in the limit  with  the corresponding $\Z^d$-versions of the Markovian and the Bosonic loop measure. The total mass of all loop measures is infinite. 

\medskip

\item[(b)] Theorem~\ref{THM-main1}:  The convergence of the finite-dimensional distributions continues to hold for the independent space-time random walk either for $ d\ge 3 $ and $ \mu\le 0 $, or for $ d\ge 1 $ and $ \mu<0 $.

\medskip

\item[(c)] Theorem~\ref{THM-main2}: The results  hold for the independent space-time random walk either for $ d\ge 3 $ and $ \mu\le 0 $, or for $ d\ge 1 $ and $ \mu<0 $ when $ F\colon[0,\infty)^{\Z^d}\to\R $ is bounded and continuous with compact support ($\supp(F)\Subset\Z^d $).

\medskip 

\item[(d)] Proposition~\ref{diagonalization_theorem}: All statements continue to hold for the independent space-time random walk in the thermodynamic limit.

\medskip

\item[(e)] Theorem~\ref{Isomorphism occupation field}: The statements continue to hold in the thermodynamic limit when either $ d\ge 1 $ and  $\mu<0$, or when  $d\ge 3$ and $ \mu\le 0 $. 

\medskip

\end{enumerate}
\end{theorem}

%
%

\section{The physical background of the Bosonic loop measure}\label{S:Bosons}
We outline the physical background of the Bosonic loop measure.  We fix our graph $ \gb=\L\Subset\Z^d $ to be a connected finite subset of $ \Z^d $ and consider the simple random walk in $ \L $ with Dirichlet boundary conditions, that is,  with killing upon hitting the boundary of $ \L $, see Section~\ref{S:intro}.  In Section~\ref{sec-bosons}   we outline the connection between  the Bosonic loop measure and the grand canonical partition function for systems of Bosons at thermodynamic equilibrium, whereas in Section~\ref{sec-BEC} we discuss the onset of Bose-Einstein condensation as  a quantum phase transition.

\subsection{Systems of Bosons in thermodynamic equilibrium}\label{sec-bosons}
In quantum mechanics particles can either be Bosons or Fermions. We will now describe the distinguishing feature of Bosons.  We consider a system of interacting Bosons on the lattice $ \L $. That is, the one-particle Hilbert space is $ \Hscr_\L=\ell_2(\L) $, and the $N$-particle Hilbert space is just the tensor product. The energy is given by the Hamilton operator for $N$ particles 
$$
H_N=-\sum_{i=1}^N \Delta_i^{\ssup{\L}}+\sum_{1\le i<i\le N} v(|x^{\ssup{i}}-x^{\ssup{j}}|)\, ,
$$
where $ \Delta_i^{\ssup{\L}} $ is the discrete Laplacian operator in $\L $ with Dirichlet boundary conditions giving the kinetic energy for particle $i$, and the interaction potential function $ v\colon\R\to\R  $  depends on the mutual lattice distance of particle $i$ at $ x^{\ssup{i}}\in\L $ and particle $j$ at $x^{\ssup{j}}\in\L $. If the particle number is not known exactly, the thermodynamic equilibrium is described  by the grand canonical ensemble where  the Hilbert space is the direct sum  
$$
\Fcal=\bigoplus_{N=0}^\infty \Hscr_\L^{\otimes N},
$$
called Fock space. 
States of identical and indistinguishable Bosons are given by symmetric wave functions, that is, for any number $ N$ of Bosons possible states  are given by all symmetric wave functions in the tensor product $ \Hscr_\L^{\otimes N} $. This symmetry  is the unique distinguishing feature of Bosons contrary to Fermions whose states are given by anti-symmetric wave functions. We write $ \Fcal_+ $ for the Fock space of all symmetric wave functions. The thermodynamic equilibrium of Bosons at inverse temperature $ \beta $ and chemical potential $ \mu $ is given  by the  grand canonical partition function which is the trace over the symmetrised Fock space, 
\begin{equation}\label{partition}
Z_{\L,v}(\beta,\mu)=\tr_{\Fcal_+}(\ex^{-\beta(H-\mu N)}),
\end{equation}
where $ H $ is the Hamilton operator having projection $ H_N $ on the sector $ \Hscr_\L^{\otimes N} $ and where $ N $ is the number operator in $ \L $ taking the value $N$ on the sector $ \Hscr_\L^{\otimes N} $.  It is the trace operation on the right hand side of \eqref{partition} which turns the problem of calculating the partition function into a classical probability question. Using the Feynman-Kac formula for traces,  we derive the following representation of the grand canonical partition function,
\begin{equation}\label{partition1}
Z_{\L,v}(\beta,\mu)=\sum_{N=0}^\infty\frac{\ex^{\beta\mu N}}{N!}\sum_{x_1\in\L}\cdots\sum_{x_N\in\L}\sum_{\sigma\in\Sym_N}\bigotimes_{i=1}^N \P_{x_i,x_{\sigma(i)}}^{\ssup{\beta}}\big[\ex^{\sum_{1\le i,j\le N}\int_0^\beta v(|X^{\ssup{i}}_t-X^{\ssup{j}}_t|)\,\d t}\big],
\end{equation}
where $ \Sym_N $ is the set of all permutations of $N$ elements, and the right hand side can be interpreted as a system of $N$ random walks $ (X^{\ssup{i}}_t)_{t\ge 0}, i=1,\ldots, N $, under symmetrised initial and terminal conditions (see \cite{AD08} for details).  Following \cite{Gin71,BR97} and \cite{ACK} we employ a cycle-expansion using the definition of  the Bosonic loop measure to obtain 
\begin{equation}\label{partition2}
Z_{\L,v}(\beta,\mu)=\sum_{N=0}^\infty\frac{1}{N!}\bigotimes_{i=1}^N M_{\L,\mu,\beta}^B(\d\omega^{\ssup{i}})\big[\ex^{-V(\omega^{\ssup{1}},\ldots,\omega^{\ssup{N}})}\big],
\end{equation}
where the interaction energy of  $N$ loops is the functional
\begin{equation}\label{interaction}
\begin{aligned}
&V(\omega^{\ssup{1}},\ldots,\omega^{\ssup{N}})=\\
&\frac{1}{2}\sum_{1\le i,j\le N}\sum_{k=0}^{[\ell(\omega^{\ssup{i}})-1]/\beta}\,\,\,\sum_{m=0}^{[\ell(\omega^{\ssup{j}})-1]/\beta}\1\{(i,k)\not=(j,m)\}\int_0^\beta\,v(|\omega^{\ssup{i}}(k\beta+t)-\omega^{\ssup{j}}(m\beta+t)|)\,\d t,
\end{aligned}
\end{equation}
and where we write $ \ell(\omega^{\ssup{i}}) $ for the length of the $i$-th loop. The interaction functional \eqref{interaction}  comprises the interactions of all loop legs (a loop $\omega$ of length $k\beta$ consists of $k$ legs, each of which is a path $(\omega(j\beta +t))_{t\in [0,\beta]}$ for $0\le j<k$) with  time horizon $ [0,\beta] $ in all Bosonic loops.
For $ v\equiv 0 $, using \eqref{BosLoopMDef} the expression reduces to the following formula
\begin{equation}\label{partition3}
Z_\L(\beta,\mu)=\exp\bigg(\sum_{x\in\L}\sum_{j=1}^\infty\frac{\ex^{\beta\mu j}}{j} \P^{\ssup{j\beta}}_{x,x}(\Gamma)\bigg)=\exp(M_{\L,\mu,\beta}^B(\Gamma))\, .
\end{equation}
This shows that the total mass of the Bosonic loop measure equals the logarithm of the grand canonical partition function for a system of non-interacting identical Bosons at thermodynamic equilibrium.
The connection goes further as we can represent (quantum) correlation functions as well. 
Quantum correlations are given as  reduced traces of the equilibrium state. 
We focus solely on the first correlation function, called the one-particle reduced density matrix $ \rho_{\L,\Phi}^{\ssup{1}} $, which is given by the partial trace after integrating out all but one of the particles, 
\begin{equation}\label{rdm1}
\rho_{\L,v}^{\ssup{1}}=\frac{1}{Z_\L(\beta,\mu)}\sum_{n=0}^\infty \tr_{\Hscr_{\L,+}^{\otimes (n-1)}}\,\Big(\ex^{-\beta( H_n-\mu n)}\Big).
\end{equation}

With \eqref{partition2} and \cite{BR97} we obtain a representation for the kernel of  the trace class operator $ \rho^{\ssup{1}}_{\L,v} $, 
\begin{equation}\label{reduced1}
\begin{aligned}
\rho^{\ssup{1}}_{\L,v}(x,y)&=\frac{1}{Z_{\L,v}(\beta,\mu)}\sum_{j=1}^\infty\sum_{N=0}^\infty \frac{\ex^{\beta\mu j}}{N!} \P_{x,y}^{\ssup{j\beta}}\otimes\big(M_{\L,\mu,\beta}^B\big)^{\otimes N}
\Big(\ex^{-V(\omega,\omega^{\ssup{1}},\ldots,\omega^{\ssup{N}})}\Big),\quad x,y\in\L.
\end{aligned}
\end{equation}
The term under the sum over $N$ is an expectation with respect to $N$ independent loops $ \omega^{\ssup{1}},\ldots, \omega^{\ssup{N}} $ chosen according to the Bosonic loop measure $ M_{\L,\beta,\mu}^B $ \eqref{BosLoopMDef}, whilst the expectation is with respect to a single random walk path from $x$ to $y $, whose length (time horizon) is determined according to the weighted sum $ \sum_{j\ge 1}\ex^{\beta\mu j}\P_{x,y}^{\ssup{ j\beta}} $. The sum over $ N\ge  0 $ is none other than the Poisson point process with intensity measure  $ M_{\L,\mu,\beta}^B $.  
Henceforth, we obtain a representation via an expectation of the Poisson loop process with intensity measure $ M^B_{\L,\mu,\beta} $.   Recall that for any counting measure $ \eta\in\Omega $ we write the interaction functional  as
$$
V(\eta)=V(\omega^{\ssup{1}},\ldots,\omega^{\ssup{N}}) \;\mbox{ for } \eta=\sum_{i=1}^N\delta_{\omega^{\ssup{i}}}.
$$

\begin{theorem}\label{physical quantities representation}
The grand-canonical partition function  and the kernel of the one-particle reduced density matrix are  given by the following functionals of the Bosonic loop process $ \Pr_{\L,\mu,\beta}^B $
\begin{equation}
\begin{aligned}
Z_{\Lambda,v}(\beta,\mu)&=\ex^{M_{\Lambda,\mu,\beta}^B(\Gamma)}\Er_{\Lambda,\mu,\beta}^B\Big[\ex^{-V(\eta)}\Big]\, =  Z_{\Lambda}(\beta,\mu)\Er_{\Lambda,\mu,\beta}^B\Big[\ex^{-V(\eta)}\Big]\,,  \\[1.5ex]
\rho^{\ssup{1}}_{\Lambda,v}(x,y)&=\frac{\E_{x,y,B}^{\mu,\beta}\otimes\Er_{\Lambda,\mu,\beta}^B\Big[\ex^{-V(X,\eta)}\Big]}{\Er_{\Lambda,\mu,\beta}^B\Big[\ex^{-V(\eta)}\Big]}\, \quad x,y\in\L,
\end{aligned}
\end{equation}
where $ \E_{x,y,B}^{\mu,\beta} $ is the expectation with respect to the measure in \eqref{bosonbridgemeasure}.
\end{theorem}

\begin{proofsect}{Proof} 
The proof follows easily from \eqref{partition2}, \eqref{partition3} and \eqref{reduced1} by noting that
$$
\Er_{\L,\mu,\beta}\big[\ex^{-V(\eta)}\big]=\sum_{N=0}^\infty\ex^{-M_{\L,\mu,\beta}^B(\Gamma)}\frac{1}{N!} \big(M_{\L,\mu,\beta}^B\big)^{\otimes N}
\Big(\ex^{-V(\omega^{\ssup{1}},\ldots,\omega^{\ssup{N}})}\Big).$$

\qed
\end{proofsect}

In the next lemma we show that the spatial mean of the occupation time in $ \L $ is in fact the partial derivative of the partition function with respect to the chemical potential, and thus, it is the expected particle density which we denote by $ \rho_\L(\beta,\mu) $.

\begin{lemma}    
$$
\Er_{\L,\mu,\beta}^B\big[\frac{1}{|\L|}\sum_{x\in\L}\Lscr_x\big]=\frac{1}{\beta}\partial_\mu\log Z_\L(\beta,\mu).
$$
\end{lemma}

\begin{proofsect}{Proof}
Recall that $ \Er_{\L,\mu,\beta}^B $ is the expectation with respect to the  Bosonic Poisson process. Then,
$$
\begin{aligned}
\Er_{\L,\mu,\beta}^B\big[\frac{1}{|\L|}\sum_{x\in\L}\Lscr_x\big]&=\frac{1}{|\L|}\sum_{x\in\L}\int\int_0^{\ell(\omega)}\,\1\{X_s(\omega)=x\}\,\d s \,M_{\L,\beta,\mu}(\d\omega)
=\frac{\beta}{|\L|}\sum_{j=1}^\infty jM_{\L,\beta,\mu}[\Gamma_{j\beta}]\\[1.5ex] &=\frac{1}{|\L|}\partial_\mu M_{\L,\beta,\mu}[\Gamma]=\frac{1}{\beta}\partial_\mu\log Z_\L(\beta,\mu)=:\rho_\L(\beta,\mu),
\end{aligned}
$$
 where $ \rho_\L(\beta,\mu) $ is the expected particle density in $ \L$. 
 \qed
 \end{proofsect}

 
\subsection{Bose-Einstein condensation - relevance and discussion}\label{sec-BEC}

We outline the original calculations by Bose and Einstein 1925 in our occupation field setting.  Detailed proofs of all  the statements below can be found in one  \cite{BR97} or \cite{LSSY05}. 
S\H{u}t\H{o}  \cite{S02} has derived a criterion for the onset of Bose-Einstein condensation in terms of losing probability mass on the distribution of loops of any finite lengths. We are not repeating his arguments here but consider the following approach. The first observation concerns the thermodynamic limit 
$$
\lim_{\L\uparrow\Z^d} \rho_\L(\beta,\mu)=\rho(\beta,\mu)\in (0,\infty),
$$
which holds in the ideal Bose gas (\cite{S02, BR97, LSSY05}) for $ \beta> 0 $ and $ \mu< 0 $. For given $\beta>0 $ we look for solutions  such that the following holds.
$$
\beta>0\colon \mbox{for } \rho>0 \mbox{ find } \mu^*=\mu^*(\rho) \mbox{ with } \rho(\beta,\mu^*)=\rho.
$$
It turns out that a solution exists only for 
$$ \rho<\rho_{\rm c}(\beta):=\lim_{\mu\nearrow 0}\rho(\beta,\mu)=\begin{cases} \infty &,\mbox{ if } d=1,2,\\<\infty &, \mbox{ for } d\ge 3.
\end{cases}
$$
For dimensions $ d\ge 3 $ and any particle density $ \rho>\rho_{\rm c}(\beta) $, the excess particle mass density $ \rho-\rho_{\rm c}(\beta) $  is interpreted as the density of the {\em Bose-Einstein condensate}.  One of the most prominent open problems in mathematical physics is the understanding of  {\em Bose-Einstein condensation (BEC)} for interacting Bosons.  This phase transition is characterised by the fact that  a macroscopic part of the system condenses to a state which is highly correlated.  
Only partial successes in proving BEC have been achieved, like the description of the free energy of the non-interacting, system (already contained  in Bose's and Einstein's seminal paper in 1925) or the analysis of mean-field models (e.g. \cite{AD18}) or  the analysis of dilute systems at vanishing temperature \cite{LSSY05} or the proof of BEC in lattice systems with half-filling \cite{LSSY05}.  In \cite{ACK}  the authors provide  a formula for the limiting free energy in  the so-called canonical ensemble where the particle density is fixed.  It turns out that the  formula is difficult to analyse, and it is one of the main aims of the current work to provide an alternative approach in terms of space-time Green functions.

 A definition for BEC for interacting gases was first provided by Onsager and Penrose \cite{OP56}, who studied quantum correlations given by the 1-particle reduced density matrix. The phase transition is signalled by the occurrence of off-diagonal long range order, meaning that the kernel $ \rho^{\ssup{1}}_{V} (x,y)=\lim_{\Lambda\uparrow\Z^d} \rho^{\ssup{1}}_{\L,V} (x,y) $ of the one-particle reduced density matrix in the thermodynamic limit $ \L\uparrow \Z^d $ has non-vanishing off-diagonal entries
 $$
 \rho^{\ssup{1}}_v(x,y)\to c>0 \;\mbox{ as } \; |x-y|\to\infty.
$$
This has been proved for the ideal Bose gas ($v\equiv 0$). The representation of the one-particle reduced density kernel in Theorem~\ref{physical quantities representation} for $ v\equiv 0 $ allows to perform the thermodynamic limit (cf. \cite{BR97}),
\begin{equation}\label{RDM}
 \rho^{\ssup{1}}(x,y)=\lim_{\L\uparrow\Z^d}\rho^{\ssup{1}}_\L(x,y)=\sum_{j=1}^\infty \ex^{\beta\mu j}\P_{x,y}^{\ssup{j\beta}}(D).
\end{equation} 

The probabilistic loop expansion of $ \rho^{\ssup{1}}(x,y) $ allows a justification of the off-diagonal long range order only with an additional assumption. Namely, for non-interacting systems one has to add the information that for large particle densities the excess particle mass is carried in cycles of ``infinite'' length. In the future we aim to identify these  ``infinitely'' long cycles in our space-time setting, either as unbounded winding numbers around the ``time'' torus, or as Bosonic interlacements. Recent progress has been made in \cite{AD18} on various mean-field type model.

\section{Proofs}\label{section proofs}
This section contains all  remaining proofs. In Subsection~\ref{Proof1} we collect all proofs for Section~\ref{Markovian and Bosonic loop measures} on all properties of our loop measures and their space time versions and torus limits, whereas  Subsection~\ref{Proof2} provides all proofs for the isomorphism theorems and moment formulae.

\subsection{Loop measure properties}\label{Proof1}

\begin{proofsect}{Proof of Proposition~\ref{properties}}
For any $k\in\N $, let $ 0\le t_1<t_2<\cdots< t_k $ and let $ A\subset (0,\infty) $ be measurable such that $ \inf(A)>t_k $. We need  $\inf(A)>t_k$ because we cannot specify distributions at  times which are not defined (larger than the lengths of each loop with $ \ell\in A$).  Thus 
\begin{equation}
\begin{aligned}
M_{\gb,\mu,\beta}^B & [X_{t_1}=x_1,\ldots,X_{t_k}=x_k,\ell\in A]=\sum_{x\in\gb}\sum_{\beta j\in A}\frac{\ex^{\beta \mu j}}{j}\P_x\left(X_{t_1}=x_1,\ldots,X_{t_k}=x_k,X_{\beta j}=x\right)\\[1ex]
&=\sum_{\beta j\in A}\sum_{x\in\gb}\frac{\ex^{\beta \mu j}}{j}p_{t_1}(x,x_1)p_{t_2-t_1}(x_1,x_2)\cdots p_{\beta j-t_k}(x_k,x)\\[1ex] &=\sum_{\beta j\in A}\frac{\ex^{\beta \mu j}}{j}p_{t_2-t_1}(x_1,x_2)\cdots p_{\beta j-t_k+t_{1}}(x_k,x_1).
\end{aligned}
\end{equation}
We have used the definition of the loop measure in the first, the Markov property in the second and the Chapman-Kolmogorov equations in the last step. The proof of the  result for the Markovian loop measure follows similarly, see details in \cite{LeJan}. This finishes the proof of \eqref{finite_Boson _distr} and \eqref{finite_Boson _distr k=1}, and thus \eqref{FDD-boson}. \\
To calculate the mass of all loops with at least one jump, i.e., to prove \eqref{mass_markov}, we note that 
\begin{equation}
\begin{aligned}
&M_{\gb,\mu}[\{ \text{ there is at least one jump } \} ]=\sum_{x\in\gb}\int_0^\infty\frac{\ex^{\mu t}}{t}\Big(\ex^{tQ}-\ex^{tD}\Big)(x,x)\d t\\[1ex]
&=\sum_{\lambda\in\gb^*}\int_0^\infty\frac{\ex^{\mu t}}{t}\Big(\ex^{tq(\lambda)}-\ex^{td(\lambda)}\Big)\d t=\sum_{\lambda\in\gb^*}\log\Big(\frac{d(\lambda)+\mu}{q(\lambda)+\mu}\Big)=\log\Big(\frac{\det(D+\mu I)}{\det(Q+\mu I)}\Big)\, ,
\end{aligned}
\end{equation}
where we used 
\begin{equation}\label{logformula}
\log\Big(\frac{a}{b}\Big)=\int_0^\infty\,\frac{\d t}{t}\Big(\ex^{-tb}-\ex^{-at}\Big)\,,\quad a>b>0.
\end{equation}
Furthermore, we used the fact that the set of generalised eigenvalues with multiplicity has the same cardinality as $ \gb $, and thus we take the sum over $ \lambda\in\gb^* $ with $ \# \gb=\# \gb^* $. In the following we denote $ q(x) $ the $x$-th generalised eigenvalue of the generator $Q$ and $ d(x) $ the $x$-th eigenvalue of the diagonal part $D$ of the matrix $ Q$.
This notation is short but slightly imprecise: There is no natural ordering or labelling of the (extended) eigenvalues of non-diagonal matrices. The sum is thus to be understood as choosing an arbitrary ordering of the eigenvalues first.  

We are left to show \eqref{mass_boson} and \eqref{mass_boson_total} for the Bosonic loop measure. For this Bosonic case we proceed similarly  as done above for the Markovian loop measure,
\begin{equation}
\begin{split}
&M_{\gb,\mu,\beta}^B[\{ \text{ there is at least one jump } \} ]=\sum_{x\in\gb}\sum_{j\ge 1}\frac{\ex^{\mu \beta j}}{j}\Big(\ex^{\beta jQ}-\ex^{\beta jD}\Big)(x,x)\\[1.5ex]
&=\sum_{x\in\gb}\sum_{j\ge 1}\frac{\ex^{\mu \beta j}}{j}\Big(\ex^{\beta jq(x)}-\ex^{\beta jd(x)}\Big)=\sum_{x\in\gb}\log\Big(\frac{1-\ex^{\beta(d(x)+\mu)}}{1-\ex^{\beta(q(x)+\mu)}}\Big)=\log\Big(\frac{\det(I-\ex^{\beta(D+\mu I)})}{\det(I-\ex^{\beta(Q+\mu I)})}\Big)\, .
\end{split}
\end{equation}
This concludes the proof of \eqref{mass_boson}. 
Calculating the mass is easier. We expand
\begin{equation}
\begin{split}
M_{\gb,\mu,\beta}^B[\Gamma]=\sum_{x\in\gb}\sum_{j\ge 1}\frac{\ex^{\mu \beta j}}{j}\ex^{\beta j Q}(x,x)=-\text{Tr}\log\Big(I-\ex^{\beta(Q+\mu I)}\Big)=-\log\Big(\det\Big(I-\ex^{\beta(Q+\mu I)}\Big)\Big)\, .
\end{split}
\end{equation}
We need $\kappa(x)-\mu>0$ for at least one $x\in\gb$ as otherwise $0$ is an element of the spectrum of $(Q+\mu)$ and the associated sums are infinite. We have shown that both sides of the equation are equal and thus verified equation \eqref{mass_boson_total}. The claim for the mass of the Markovian loop measure is proved in \cite{LeJan}, page 21.
\qed
\end{proofsect}

\medskip

\begin{proofsect}{Proof of Theorem~\ref{THM-main1}}
The proof is organised in three steps. In the first step we show the convergence \eqref{Loop_measure_convergence_theorem} for finite dimensional distributions. In the second step we show the weak convergence, and in the third step we justify tightness which is needed in step 2.

\medskip

\noindent \textbf{Step 1:} \\[1ex]
We first prove the convergence \eqref{Loop_measure_convergence_theorem} for the independent space-time random walk. Let us first examine the case with  $k=1$ as it is different  from the cases with $k>1$ and helps to develop intuition. By using the independence of the space and the time coordinates, the Markov property, the formula for the finite dimensional distributions in Proposition~\ref{properties} and the translation invariance on the torus $\T_N$,   we can rewrite the left hand side of \eqref{Loop_measure_convergence_theorem} as
\begin{equation}\label{loop_measure_convergence_k=1}
\begin{aligned}
M_N&[X_{t_1}=x_1,\ell\in A]=\sum_{\tau_1=0}^{N-1}M_{\Lambda\times\T_N,\mu}\left[(X_{t_1}^{\ssup{1}},X^{\ssup{2}}_{t_1})=(x_1,\tau_1),\ell\in A\right]\\[1.5ex]
&=\sum_{\tau_1=0}^{N-1}\sum_{\tau=0}^{N-1}\sum_{x\in\L}\int_0^\infty\frac{\ex^{\mu t}}{t}\,\P^{\ssup{t}}_{(x,\tau),(x,\tau)}\big((X_{t_1}^{\ssup{1}},X^{\ssup{2}}_{t_1})=(x_1,\tau_1),t\in A\big)\,\d t,\\[1.5ex]
& = \int_Ap_t(x_1,x_1)\ex^{-tN\beta^{-1}}N\Big(\sum_{j=0}^\infty\frac{(tN\beta^{-1})^{jN}}{(jN)!}\Big)\frac{\ex^{t\mu}}{t}\,\d t\\[1.5ex]& =
\beta\sum_{j=0}^\infty\int_A p_t(x_1,x_1) t^{jN} \ex^{-tN\beta^{-1} }\frac{(N\beta^{-1})^{jN+1}}{\Gamma(jN+1)}\frac{\ex^{t\mu}}{t}\d t\\[1.5ex]
&=\beta\sum_{j=0}^\infty\E_{jN+1,N\beta^{-1}}\left[\1_A(X)p_X(x_1,x_1)\frac{\ex^{X\mu}}{X}\right]\, ,
\end{aligned}
\end{equation}
since the independence of torus and spatial component implies $\P_{(x_1,\tau)}(X_t=(x_1,\tau))=p_t(x_1,x_1)\P^{\ssup{t}}_{N,\tau,\tau} $, where $ \P^{\ssup{t}}_{N,\tau,\tau} $ is the random walk bridge measure for the torus component $ X^{\ssup{2}}_t $. The expectation $\E_{jN+1,N\beta^{-1}}$ is with respect to a Gamma distributed random variable $X$ with shape parameter $(jN+1)$ and rate $N\beta^{-1}$.   For the exchange of the limit and the sum, see equation \eqref{exponentialest}. The distribution of the torus component is invariant under torus translations, implying 
$$
\sum_{\tau=0}^{N-1}\P^{\ssup{t}}_{N,\tau,\tau}=N\P^{\ssup{t}}_{N,0,0}=N\ex^{-tN\beta^{-1} }\sum_{j=0}^\infty\frac{(tN\beta^{-1})^{jN}}{(jN)!}  .
$$
Since $\kappa(x)-\mu>0$ for at least one $x\in\Lambda$, note that $p_t(x,y)\ex^{t\mu}$ decays exponentially in $t$ (in case $ \mu=0 $ the decay comes from the killing term, see \cite{Sznitman}). We use that $\inf(A)>t_1>0 $ to bound the denominator with some $C>0 $, and for some $\alpha>0$ small enough,
\begin{equation}\label{exponentialest}
\begin{aligned}
\E_{jN+1,N\beta^{-1}}\Big[\1_A(X)p_X(x_1,x_1)\frac{\ex^{X\mu}}{X}\Big]& \le C \E_{jN+1,N\beta^{-1}}\Big[\ex^{-\alpha X}\Big]=C\Big(1+\frac{\alpha\beta}{N}\Big)^{-Nj+1}\\[1ex] & \le C\left(\frac{\ex^{-\beta \alpha}}{2}\right)^j\, ,
\end{aligned}
\end{equation}
and,  as the bound is independent of $N$ and decays exponentially in $ j $, the series converges absolutely uniformly in $N$.  Thus we can approximate the series with a finite sum in  \eqref{loop_measure_convergence_k=1}, or, equivalently, can commute the summation with differentiating the expectation value for every $ j\in\N $.
For any $ j\in\N $, one then expands the integrand around $\frac{jN+1}{N\beta^{-1}}=\E_{jN+1,N\beta^{-1}}\left[X\right]$, i.e., for $g(t)=\1_A(t)p_t(x_1,x_1)\frac{\ex^{t\mu}}{t}$, which is differentiable infinitely often (for $N$ large as $ \inf(A)>0 $),  one gets
\begin{equation}
\begin{split}
&\E_{jN+1,N\beta^{-1}}\Big[\mathbbm{1}_A(X)p_X(x_1,x_1)\frac{\ex^{X\mu}}{X}\Big]\\[1.5ex]
&=g\Big(\frac{jN+1}{N\beta^{-1}}\Big)+g^\prime\Big(\frac{jN+1}{N\beta^{-1}}\Big)\E_{jN+1,N\beta^{-1}}\Big[\Big(X-\frac{jN+1}{N\beta^{-1}}\Big)\Big]+\frac{jN+1}{(N\beta^{-1})^2}\,\Ocal(1)\, .
\end{split}
\end{equation}
The variance is  $\frac{jN+1}{(N\beta^{-1})^2}$ and goes to zero for $N\to\infty$. Higher moments of the Gamma distribution are of the order $ \frac{1}{N}\Ocal(1) $. Since the first moment vanishes and the point of expansion converges to $j\beta$ uniformly for all $ j\in\N $ as  $N\to\infty$, we get (the $j=0 $ term vanishes as $ \inf(A)\ge t_1>0 $)
\begin{equation}
\lim_{N\to\infty} M_{N}\left[X_{t_1}=x_1, \,\ell\in A\right]  =  \beta\sum_{j=1}^\infty p_{\beta j}(x_1,x_1)\mathbbm{1}_{A}(\beta j)\frac{\ex^{\beta j \mu}}{\beta j}=\sum_{j=1}^\infty p_{\beta j}(x_1,x_1)\frac{\ex^{\beta\mu j}}{j}\1_A(\beta j)\, .
\end{equation}
Note that the requirement $\beta\N\cap\partial  A=\emptyset$ was necessary for us to expand around $\frac{jN+1}{N\beta^{-1}}$ (for $N$ large enough) as otherwise the function $g(t)$ has a jump at $j\beta, j\in\N$, and thus would not be differentiable.

\medskip

We turn to the case $k\ge 2$ and introduce some notation. Denote $ p^{\ssup{k}} $ the product of the jump probabilities,
\begin{equation}
p^{\ssup{k}}=\prod_{j=1}^{k-1}p_{t_{j+1}-t_j}(x_j,x_{j+1})\, .
\end{equation}
Write the left hand side of \eqref{Loop_measure_convergence_theorem} using Chapman Kolmogorov again as 
\begin{equation}\label{long_equation_boson_bproof}
\begin{aligned}
M_N&[X_{t_1}=x_1,\ldots, X_{t_k}=x_k]=\sum_{\tau=0}^{N-1}\sum_{x\in\L} \int_0^\infty\,\frac{\ex^{\mu t}}{t}\,\P^{\ssup{t}}_{(x,\tau),(x,\tau)}\big(\pi_\L(X)_{t_1}=x_1,\ldots, \pi_\L(X)_{t_k}=x_k\big)\\[1.5ex]
&=p^{\ssup{k}}\sum_{\heap{\tau_1,\ldots,\tau_k}{\in\{0,\ldots,N-1\}}}\Big(\prod_{m=1}^{k-1}\P^{\ssup{t_{m+1}-t_{m}}}_{N,\tau_{m+1},\tau_{m}}\Big)\times\int_A\,\Big(p_{t-t_k+t_1}(x_k,x_1)\P^{\ssup{t_{1}+t-t_k}}_{N,\tau_{k},\tau_{1}}\frac{\ex^{t\mu}}{t}\Big)\,\d t\\[1.5ex]
&=p^{\ssup{k}} \int_A \,\Big( p_{t-t_k+t_1}(x_k,x_1)N\P^{\ssup{t}}_{N,0,0}\frac{\ex^{t\mu}}{t}\Big)\,\d t\\
\end{aligned}
\end{equation}

As in the case $k=1 $ we obtain absolute uniform convergence of the series. Similar to the case $k=1$ we now perform a second order Taylor expansion of $g(t)$ around $t=\frac{jN+1}{N\beta^{-1}}$. We obtain that the left hand side of \eqref{Loop_measure_convergence_theorem} converges, as $ N\to\infty $, to
\begin{equation}
p^{\ssup{k}}\beta \sum_{j=0}^\infty \mathbbm{1}_A(n\beta)p_{j\beta-t_k+t_1}(x_k,x_1)\frac{\ex^{\mu j\beta}}{j\beta}\, .
\end{equation}
This finishes the proof of \eqref{Loop_measure_convergence_theorem} in Theorem~\ref{THM-main1} for independent space-time random walks who are strongly asymptotically independent space-time random walks.

We now turn to proof of statement \eqref{Loop_measure_convergence_theorem} for general  strongly asymptotically independent space-time random walks with generator $G_N=Q_N+E_N$ and probability $ \widetilde{\P}_{(x,\tau)} $. Denote $\widetilde{M}_N$ the projected space-time loop measure induced by $G_N$. We shall estimate the difference between $ \widetilde{M}_N $ and  $M_N $, where $ M_N$ is  the projected loop measure for the independent space-time random walk with generator $Q_N $.  Similar to the previous part of the proof, we can find for any $ \epsilon>0 $ a parameter $M>0$ such that, independent of $N$, we obtain  
$$
M_N[\ell>M]+\widetilde{M}_N[\ell>M]<\epsilon.
$$ 

Abbreviate for $m=1,\ldots,k-1$, $x_m,x_{m+1}\in\Lambda$, $\tau_m,\tau_{m+1}\in\T_N$ and $t_m<t_{m+1}$,
\begin{equation}
\begin{aligned}
p(m)&:=\ex^{(t_{m+1}-t_m)Q_N}(x_m,\tau_m,x_{m+1},\tau_{m+1}) \quad\text{ and }\\[1.5ex]
g(m):&=\ex^{(t_{m+1}-t_m)G_N}(x_m,\tau_m,x_{m+1},\tau_{m+1})\, .
\end{aligned}
\end{equation}
For $m=k$ define
\begin{equation}
p(k):=\int_A\frac{\ex^{t\mu}}{t}\ex^{(t-t_{k}+t_1)Q_N}(x_k,\tau_k,x_{1},\tau_{1})\d t
\end{equation}
and $g(k)$ analogously. We estimate, setting $\tau_{k+1}:=\tau_1$, 
\begin{equation}\label{long calculation class one}
\begin{aligned}
&|\widetilde{M}_N\left[X_{t_1}=x_1,\ldots,X_{t_k}=x_k,\ell\in A\right]-M_N\left[X_{t_1}=x_1,\ldots,X_{t_k}=x_k,\ell\in A\right]|\\[1.5ex]
&=\Big|\sum_{\heap{\tau_1,\ldots,\tau_k}{\in \{0,\ldots,N-1\}}}\Big(\prod_{j=1}^kp(j)-\prod_{j=1}^kg(j)\Big)\Big|\\[1ex]
&=\Big|\sum_{\heap{\tau_1,\ldots,\tau_k}{\in \{0,\ldots,N-1\}}}\sum_{m=1}^k\Big(\prod_{j=1}^{m-1}p(j)\Big)(g(m)-p(m))\Big(\prod_{j=m+1}^kg(j)\Big)\Big|\\[1.5ex]
&=\sum_{m=1}^k\sum_{\tau_m,\tau_{m+1}=1}^{N-1}\big|g(m)-p(m)\big|\sum_{\heap{\tau_1,\ldots,\tau_k\setminus \{\tau_m,\tau_{m+1}\}}{\in \{0,\ldots,N-1\}}}\Big(\prod_{j=1}^{m-1}p(j)\Big)\Big(\prod_{j=m+1}^kp(j)\Big)\\[1.5ex]
&\le\sum_{m=1}^{k-1}\sum_{\tau_m,\tau_{m+1}=1}^{N-1}\big|(p(m)-g(m))\big|  \sum_{\tau=0}^{N-1}\P_{(x_1,\tau)}(X_{t_m}=(x_m,\tau_m))  \Big\{  \\[1.5ex]
&\quad \int_{A\cap (0,M]}\,\d t \;\Big(\frac{\ex^{t\mu}}{t}\,  \widetilde{\P}_{(x_{m+1},\tau_{m+1})}(X_{t-t_m+t_1}=(x_1,\tau))\Big)\Big\}   +R_k\\[1.5ex]
&\le C \sum_{m=1}^{k-1}\norm{\ex^{(t_{m+1}-t_m)Q_N}(x_m,x_{m+1})-\ex^{(t_{m+1}-t_m)G_N}(x_m,x_{m+1})}_1\times \\[1.5ex] & \quad\times \max_{\tau=0,\ldots, N-1}{\P}_{(x_1,\tau)}(X_{t_m}=(x_m,\tau_m)) +R_k\, ,
\end{aligned}
\end{equation} 
where the integral has been absorbed in the constant $ C $, and where, by the Chapman-Kolmogorov equations the error term $R_k$ is given by
\begin{equation}\label{errorXZ}
    R_k=\int_{A\cap (0,M]}\frac{\ex^{t\mu}}{t}\norm{\ex^{(t-t_{k}+t_1)Q_N}(x_k,x_{1})-\ex^{(t-t_{k}+t_1)G_N}(x_k,x_{1})}_1\d t\, .
\end{equation}
Here, $\widetilde{\P}_{x_1,\tau}  $ is the distribution of  the strongly asymptotically independent space-time random walk with initial state $ (x_1,\tau) $. We have also used the notation $\ex^{tQ}(x,y)$ for the  $ \T_N\times\T_N $ matrix  $\ex^{tQ}(x,y)(\tau,\sigma)=\ex^{tQ}(x,\tau,y,\sigma)$ for $\tau,\sigma\in\T_N$.
To bound the difference of the two matrix exponentials in \eqref{long calculation class one} and \eqref{errorXZ}, we use the basic inequality
\begin{equation}
\norm{\ex^{X}-\ex^{X+Y}}_1\le \norm{Y}_1\ex^{\norm{X}_1}\ex^{\norm{Y}_1}\, ,
\end{equation}
for two matrices $X$ and $Y$. In fact, this  inequality holds for every sub-multiplicative matrix norm. A bound for norm of the $ \T_N\times\T_N $ matrix $Q_N(x,y) $ is given by
\begin{equation}
\norm{Q_N(x,y)}_1\le \lambda_x N(N\beta^{-1})\,,\quad x,y\in\L .
\end{equation}
Combining the above bounds with equation \eqref{long calculation class one} and \eqref{class1} we conclude with the statement \eqref{Loop_measure_convergence_theorem} for strongly asymptotically independent space-time random walks. This finishes Step 1 of our proof. 

\medskip

\noindent\textbf{Step 2:} \\[1ex]
We turn to the proof of the weak convergence of the projected loop measure towards the Bosonic loop measure in the limit $ N\to \infty $. The strategy here is to consider loops whose length is bounded away from zero such that the loop measure becomes a finite measure (see \cite{Sznitman}) as our graph is finite and the integral is bounded as well, e.g., for $ \gamma>0 $,
$$
M_N[\ell>\gamma]\le \sum_{x\in\Lambda}\int_\gamma^\infty \P_{x,x}^t(\Gamma) N\P^{\ssup{t}}_{N,0,0}\,\d t\underset{N\to\infty}{\longrightarrow} \sum_{j\ge 1\colon j\beta>\gamma}\frac{\ex^{\beta j \mu}}{j}\P_{x,x}^{\beta j}(\Gamma)<\infty.
$$
We conclude with the weak convergence once we combine Step 1, the above reasoning and the tightness to show below in Step 3.
\medskip

\noindent\textbf{Step 3:}\\[1ex]
We are left to show tightness for the sequences $M_N$ and $\widetilde{M_N}$. To start with $M_N$, fix a sequence $ (r_m)_{m\in\N} $ of positive numbers with $r_m\uparrow +\infty$ as $ m\to\infty $ with  $r_m\neq \beta n$ for all $n\in\N$. According to \cite{Bl68} it suffices to show that for all $\epsilon> 0 $ and all $m\in \N$ we have
\begin{equation}\label{tightnessC}
\lim_{\delta\downarrow 0}\limsup_{N\to\infty}M_N\left[ d'_\omega(m,\delta)\ge \epsilon\right] = 0\, ,
\end{equation}
where 
\begin{equation}
d'_\omega(m,\delta)=\inf_{\{t_i\}}\max_{1\le i\le \nu}\sup_{t,s\in [t_i,t_{i-1})}|\omega(t)-\omega(s)|\, ,
\end{equation}
and where $\{t_i\}=\big\{0=t_0<\cdots<t_\nu=r_m; |t_i-t_{i-1}|>\delta>0; \nu\in\N\big\}$ is the collection of all ordered finite subsets of $ (0,\infty) $.
The spatial jump rates $\lambda_x$ are finite, bounded from below and above, and thus \eqref{tightnessC}  holds for the projected loop measure $M_N$ of the independent space-time random walk as  the 
jumps on the torus are disregarded because of the projection. For $\widetilde{M}_N$ we note that there are three types of jumps on the space-time graph: jumps  $(x,\tau)\to(x,\tau+1)$ with rate  $N\beta^{-1}+o(\ex^{-N^2})$, jumps  $(x,\tau)\to(y,\tau)$ with rate $O(1)+o(\ex^{-N^2})$ and all other jumps with rate $o(\ex^{-N^2})$. The first type of jumps can be neglected for our purpose as they cannot be detected under the projection. For all the remaining possible  jumps, we can  estimate their jump rates  by $|\Lambda|O(1)+|\Lambda|(N+1)o(\ex^{-N^2}) $. Therefore, all rates remain bounded for $\widetilde{M}_N $ as well,  and thus  the tightness continues to hold in that case. 
\qed
\end{proofsect}

\begin{proofsect}{Proof of Theorem~\ref{THM-main2}}
As in the previous proof, we begin by showing  our statements (a)--(c)  for the independent space-time random walk first.
Using the independence of the spatial coordinate $ X^{\ssup{1}}_t $ and the torus coordinate $ X^{\ssup{2}}_t $, it is straightforward to verify that
\begin{equation}\label{splitting}
\begin{split}
M_N[F(L)]&=\sum_{x\in\L}\sum_{\tau=0}^{N-1}\int_0^\infty\frac{\ex^{t\mu}}{t}\,\E^{\ssup{t}}_{(x,\tau)}\Big[F\big(\big(\sum_{\sigma=0}^{N-1}\Lb_{x,\sigma}\big)_{x\in\Lambda}\big)\mathbbm{1}\{(X_t^{\ssup{1}},X^{\ssup{2}}_t)=(x,\tau)\}\Big]\d t\\[1.5ex]
&=\sum_{x\in\L}\sum_{j=0}^\infty\int_0^\infty\frac{\ex^{t\mu}}{t}N\E_{(x,0)}^{\ssup{t}}\Big[F\big(\big(\sum_{\sigma=0}^{N-1}\Lb_{x,\sigma}\big)_{x\in\Lambda}\big)\mathbbm{1}\{(X_t^{\ssup{1}},X^{\ssup{2}}_t)=(x,0),\wind=j\}\Big]\d t\\[1.5ex]
&=\sum_{x\in\L}\sum_{j=0}^\infty\int_0^\infty\frac{\ex^{t\mu}}{t}N\E_{x,x}^{\ssup{t}}\Big[F\big(\big(L_{x}\big)_{x\in\Lambda}\big)\Big]\P^{\ssup{t}}_{N,0}(X^{\ssup{2}}_t=0,\wind=j)\d t\,, 
\end{split}
\end{equation}
where the random variable $ \wind$  is the winding number of the torus coordinate $ X^{\ssup{2}}_t $, i.e., the number of times that $ X^{\ssup{2}}_t $  surrounds the torus, and where $\P_{N,0}^{\ssup{t}} $ is the probability with respect to the torus coordinate $ X^{\ssup{2}}_t $ with start at $0\in\T_N$. The right hand side of \eqref{splitting} can be split into the term with $ j=0 $, which means that there is no torus jump, and the term with the sum over all $ j\ge 1 $. We first show that  the latter term converges to the Bosonic loop measure. 
This can be seen as follows.  The $N\P^{\ssup{t}}_{N,0}(X^{\ssup{2}}_t=0,\wind=j)=N\P^{\ssup{t}}_{N,0,0}(\wind=j) $ can be rewritten as the density of a Gamma distributed random variable $X$ with mean $\frac{jN+1}{N\beta^{-1}}$, analogous to the proof of Theorem \ref{THM-main1}. Hence, for  each $x\in\L $ and every $ j\in\N_0 $, we write the integral with respect to $t $ as
$$
\beta \E_{jN+1,N\beta^{-1}}\Big[  \frac{\ex^{\mu X}}{X}  \E_{x,x}^{\ssup{X}}\big[F\big(\big(L_{x}\big)_{x\in\Lambda}\big)\big]\Big].
$$
Since the variance of $X$ converges to zero, we have that $X\to\beta j$ as $ N\to\infty $ in $L^2$, which in turn implies that the Gamma distribution converges (in distribution) to the delta measure $ \delta_{j\beta} $. Using our arguments for the convergence of the series in the proof of Theorem~\ref{THM-main1}, we conclude for loops with lengths $ \beta j , j\in\N $,  with the statement (a) without the second term on the right hand side of \eqref{Con}.  We shall now incorporate the missing term on the right hand side of \eqref{Con}. 

For the term with $ j=0 $, for each  $x\in\Lambda$ denote $t_x=t\delta_x$, i.e., the vector with $t$ on the $x$-th position.  Also, recall that $ (d_x)_{x\in\L} $ are the diagonal elements of the spatial generator matrix $ Q $.   Then, for each $x\in\L $, we have (using the assumptions on $F$),
\begin{equation}\label{derivative term}
\begin{aligned}
&\int_0^\infty\frac{\ex^{t\mu}}{t}N\E_{x,x}^{\ssup{t}}\big[F\big(\big(L_{x}\big)_{x\in\Lambda}\big)\big]\P^{\ssup{t}}_{N,0}(X^{\ssup{2}}_t=0,\wind=0)\,\d t\\[1.5ex]
&=\beta\int_0^\infty\frac{\ex^{t\mu}}{t}(\beta^{-1}N)\ex^{-t\beta^{-1}N}\big(F(t_x)\ex^{td_x}+(1-\ex^{td_x})\E_{x,x}^{\ssup{t}}\big[F\left(L_{x}\right)_{x\in\Lambda}\mathbbm{1}\{\text{at least one jump}\}\big]\big)\,\d t\\[1.5ex]
&=\beta\int_0^\infty\ex^{t\mu}(\beta^{-1}N)\ex^{-t\beta^{-1}N}\left((\partial_xF)(0)  \ex^{td_x}+o(t)  \right)\,\d t\\[1.5ex]
&   \; \underset{N\to\infty}\longrightarrow  \; \beta (\partial_xF)(0)\, .
\end{aligned}
\end{equation}
For the last step we use the fact that $(\beta^{-1}N)\ex^{-t\beta^{-1}N}$ is the density of an exponentially distributed random variable $X_N$ with expectation $\beta/N$. As its variance is $(\beta/N)^2$ is converging to zero for  $ N\to\infty $, we get  that $X_N\to \delta_0$ as $ N\to\infty$. We conclude with statement (a) for the independent space-time random walk. We show (a) for the other cases below after showing (b) and (c) for the independent space-time random walk. \\[1ex]
(b) For the independent space-time random walk we can insert indicators in \eqref{splitting} and \eqref{derivative term} to obtain the limits \eqref{convergence split} for the splitting of the local times. In \eqref{splitting}, for  the part with $ \wind\ge 1 $, we insert an indicator for the spatial component to be genuine loop or point loop. As $ \wind\ge 1 $, we only have genuine torus jumps as to get around the torus once requires $ N $ jumps. In \eqref{derivative term} it follows easily that no torus jump but a jump on $ \L $ has no local time, in the limit. The derivative term appears for the remaining case as shown above.\\[1ex]
(c) We prove the convergence of the occupation field in law by comparing Laplace transforms and using the previous result. We have  by Campbell's formula
\begin{equation}
\Er_{N}[\ex^{-\langle v,\mathcal{L}\rangle}]=\exp\big\{-\int_\Gamma \big(1-\ex^{-\langle v,L\rangle} \big)\d M_N\big\},\;v\in[0,\infty)^\L,
\end{equation}
and using (a) with $ F(L):=1-\ex^{\langle v,L\rangle} $,
\begin{equation}
\lim_{N\to\infty} \exp\big\{-\int_\Gamma \big(1-\ex^{-\langle v,L\rangle} \big)\d M_N\big\}=\exp\big\{-\beta\langle v,1\rangle-\int_\Gamma\big(1-\ex^{-\langle v,L\rangle} \big)\d M_{\Lambda,\mu,\beta}^B\big\}=\Er_{\Lambda,\mu,\beta}^B[\ex^{-\langle v,\mathcal{L}+\beta\rangle}]\, ,
\end{equation}
we conclude with (c) for the independent space-time random walk.\\[1.5ex]
We now prove that all statements (a)--(c) hold for general weakly asymptotically independent space-time random walks with generators $ G_N =Q_N+E_N$. We begin with the occupation field $\Lscr$ to prove (c) first and compute Laplace transforms. We then rewrite for $v\in[0,\infty)^\Lambda$ (denoting $\widetilde{\Er}_{N}$ the Poisson  process with intensity  measure $\widetilde{M}_N$ given by the space-time random walk with generator $G_N $),
\begin{equation}
\begin{aligned}
\widetilde{\Er}_{N}[\ex^{-\langle v,\Lscr\rangle}]=\exp\Big\{-M_N[1-\ex^{-\langle v,L\rangle}]-\left(\widetilde{M}_N[1-\ex^{-\langle v,L\rangle}]-M_N[1-\ex^{-\langle v,L\rangle}]\right)\Big\}\, ,
\end{aligned}
\end{equation}
Let $v_N\in[0,\infty)^{\Lambda\times\T_N}$ be defined by $v_N(x,\tau)=v(x),\, x\in\L $ such that $\langle v,\pi_\Lambda(\bL)\rangle=\langle v_N,\bL\rangle$. Hence,
\begin{equation}
\begin{aligned}
\exp\Big\{-\widetilde{M}_N[1-\ex^{-\langle v,L\rangle}]\Big\}&=
\exp\Big\{-\widetilde{M}_{\Lambda\times\T_N}[1-\ex^{-\langle v_N,\Lb\rangle}]\Big\}\\[1.5ex]&=
\exp\Big\{-\sum_{x\in\Lambda}\sum_{\tau=0}^{N-1}\int_0^\infty\frac{\ex^{t\mu}}{t}\Big(\ex^{tg_N{(x,\tau)}}-\ex^{tg_N^v{(x,\tau)}}\Big)\d t\Big\}\\[1.5ex]
&=\frac{\det(G_N+V_N+\mu I)}{\det(G_N+\mu I)}\, ,
\end{aligned}
\end{equation}
where $g_N^v{(x,\tau)}$ (respectively $g_N(x,\tau)$) is the eigenvalue indexed by $(x,\tau)$ of the matrix $G_N+V_N$ (respectively, $G_N$) with  $V_N$ being  the matrix with $v_N$ on the diagonal. The last step in the equation above follows analogously to the proof of the second half of Proposition~\ref{properties}. We shall use the following well known perturbative bound (\cite{IR08} ) for matrices $A,E\in\R^{m\times m}$
\begin{equation}
|\det(A+E)-\det(A)|\le \sum_{i=1}^m s_{m-i}\norm{E}_2^i\, ,
\end{equation}
where for $\sigma_1\ge\ldots\ge \sigma_m\ge 0$ being the singular values of the matrix $A$ we define
\begin{equation}
s_k=\sum_{1\le i_1<\ldots <i_k\le m}\sigma_{i_1}\cdots\sigma_{i_k}\,  ,
\end{equation}
and $s_0=1$. We used $\norm{E}_2$ to denote  the largest singular value of $E$ or equivalently the largest eigenvalue of $\sqrt{EE^*}$.\\
One can now bound, setting $m=|\Lambda|N$,
\begin{equation}
\begin{split}
|\det(Q_N+\mu I+E_N)-\det(Q_N+\mu I)|\le \sum_{i=1}^{m}s_{m-i}\norm{E_N}_2^i\, ,
\end{split}
\end{equation}
Using the Gershgorin circle theorem \cite{Ger31}, we can estimate
\begin{equation}
s_k\le \binom{m}{k}2^k(\lambda+N\beta^{-1}-\mu)^k\, ,\quad k=0,1,\ldots, m-1,
\end{equation}
where $\lambda:=\max_{x\in\Lambda}\lambda_x$. Thus, by the Binomial theorem and $\norm{E_N}_2$ small
\begin{equation}
\begin{split}
2^{-m}|\det(Q_N+\mu+E_N)-\det(Q_N+\mu)|&\le \Big(\lambda+N\beta^{-1}-\mu+\norm{E_N}_2\Big)^m-(\lambda+N\beta^{-1}-\mu)^m\\[1.5ex]  &\le\Ocal(m) (\lambda+N\beta^{-1}-\mu)^{m-1}\norm{E_N}_2
\end{split}
\end{equation}
We thus conclude that if the norm $\norm{E_N}_2$ satisfies the following estimate for all $\alpha>0$,
\begin{equation}\label{determinant bounds, alpha}
\begin{split}
\norm{E_N}_2=o\Big(\big(\alpha+2\beta^{-1}N\big)^{-N|\Lambda|}\Big)\, ,
\end{split}
\end{equation}
the difference between the determinants converges to zero.
By the bound \eqref{determinant bounds, alpha} for all weakly asymptotically independent space-time random walks with \eqref{class2}, using the Gershgorin circle theorem, we obtain that 
\begin{equation}
|\det(Q_N+\mu I+V_N+E_N)-\det(Q_N+\mu I+V_N)|\longrightarrow  0\, \mbox{ as } N\to\infty.
\end{equation}
Note that we used the  bound $\norm{E}_2\le \norm{E}_1$. The convergence of the determinants implies the convergence of the Laplace transforms of the space-time Poisson processes to the Poisson process with intensity measure $ M^B_{\L,\mu,\beta} $. This shows (c) for the weakly asymptotically independent space-time random walk.\\[1ex]
We are now proving part $(a)$ and $(b)$ for the weakly asymptotically independent case. Using the convergence of the Laplace transforms of the Poisson process, the convergence of the local distributions under  $\widetilde{M}_N$ can be seen as follows.  As the intensity measure $ \widetilde{M}_N $ is not finite it is difficult to relate the occupation field and the local time vector, i.e., to benefit from the convergence in (c). Our strategy here is to distinguish between genuine and point loops respectively, followed by an extension to complex-valued functions  taking advantage of  analyticity properties. Let  $J $ be the number of  jumps within a loop, i.e.,
$$
J(X):=\#\{t\colon X_t\neq\lim_{s\uparrow t}X_s\}\,. 
$$
The loop measure for at least one jump is
\begin{equation}
M_{\Lambda\times\T_N,\mu}[J>0]=\sum_{x\in\Lambda}\sum_{\tau=0} ^{N-1}\,\int_0^\infty\frac{\ex^{t\mu}}{t}\P_{(x,\tau),(x,\tau)}^t\left(\{\text{ at least one jump until time }t\,\}\right)\d t\, .
\end{equation}
Let us thus first consider the distribution of $\sum_{\tau=0}^{N-1}\Lb_{(x,\tau)}$ under the finite measure
\begin{equation}
\widetilde{M}_{>0,N}:= M_{\Lambda\times\T_N,\mu}\big(\cdot\1\{J>0\}\big)\circ\pi_\Lambda^{-1}\, .
\end{equation}
Similar to the above, we can show again with \eqref{determinant bounds, alpha} that the distribution of the occupation field $\sum_{\tau=0}^{N-1}\bLscr_{(x,\tau)}$ under  the Poisson  process $ \widetilde{\Pr}_{\Lambda\times\T_N},\mu$ with intensity measure $ \widetilde{M}_{>0,N} $ converges,
\begin{equation}\label{PPpconv}
\widetilde{\Pr}_{\Lambda\times\T_N,\mu} \circ\pi_\Lambda^{-1}\circ\Lscr\to\Pr^B_{\Lambda,\mu,\beta}\circ\Lscr\, \mbox{ as } N\to\infty.
\end{equation}
If we want to deduce from the convergence of the occupation field the convergence of the local time, we need to use a trick. For $z\in\C$ consider now the ``complexified" measure Poisson point process with expectation $\widetilde{\Er}_{\Lambda\times\T_N}^{z}$ defined as
\begin{equation}\label{Complex representation}
\begin{aligned}
\widetilde{\Er}_{\Lambda\times\T_N}^{z}&\Big[F\Big(\Big(\sum_{\tau=0}^{N-1}\bLscr_{(x,\tau)}\Big)_{x\in\Lambda}\Big)\1\{J>0\}\Big]\\[1.5ex]
&:=\ex^{-z\widetilde{M}_{>0,N}[\Gamma]}\sum_{n=0}^\infty \frac{z^n}{n!}
\Big(\widetilde{M}_{>0,N}\Big)^{\otimes n}\Big[F\Big(\Big(\sum_{i=1}^n\sum_{\tau=0}^{N-1}\bLscr_{(x,\tau)}(\bomega_i)\Big)_{x\in\Lambda}\Big)\1\{J>0\} \Big]\, ,
\end{aligned}
\end{equation}
where $\bomega_i$ is a space-time loop. For all $z\in\C\cap(0,\infty)$ the convergence of $\widetilde{\Er}^{z}_{\Lambda\times\T_N}$ to $\Er^{B,z}_{\Lambda,\beta}$ has been established in equation \eqref{PPpconv}. For fixed $F$ (bounded and continuous) the above representation is analytic in $z$ and, by the triangle inequality, bounded uniformly in $N$ on compact subsets of $\C$. Thus, as we have convergence of the sum in \eqref{Complex representation}, the individual coefficients converge as well, 
\begin{equation}
\begin{aligned}
\lim_{N\to\infty}\Big(\widetilde{M}_{>0,N}\Big)^{\otimes n}& \Big[F\Big(\Big(\sum_{i=1}^n\sum_{\tau=0}^{N-1}\bLscr_{(x,\tau)}(\bomega_i)\Big)_{x\in\Lambda}\Big)\1\{J>0\} \Big]\\[1.5ex]
&=\Big(M_{\Lambda,\beta}^B\Big)^{\otimes n}\Big[F\Big(\Big(\sum_{i=1}^n\Lscr_{x}(\omega_i)\Big)_{x\in\Lambda}\Big)\Big]\, .
\end{aligned}
\end{equation}
Now what remains is to show that for the measure
\begin{equation}
\widetilde{M}_{0,N}:=M_{\Lambda\times\T_N,\mu}\big(\cdot\1\{J=0\}\big)\circ\pi_\Lambda^{-1}\, ,
\end{equation}
we have that for every $F$ differentiable at $0$ and $F(0)=0 $,
\begin{equation}
\lim_{N\to\infty}\widetilde{M}_{0,N}[F]=\beta\sum_{x\in\Lambda}\frac{\partial}{\partial x}F(0)\, .
\end{equation}
This can be shown as follows. Denote $ \widetilde{d}(x,\tau) $ the diagonal elements of the generator $ G_N $. We compute
\begin{equation}
\begin{split}
\widetilde{M}_{0,N}[F]=\sum_{x\in\Lambda}\sum_{\tau=0}^{N-1}\int_0^\infty\frac{\ex^{t\mu}}{t}\ex^{t \widetilde{d}(x,\tau))}F(t_x)\d t=\sum_{x\in\Lambda}\int_0^\infty\frac{\ex^{t\mu}}{t}F(t_x)N\ex^{-tN\beta^{-1}}\ex^{-o(1)t}\d t\, .
\end{split}
\end{equation}
We obtain the statement by the similar  reasoning as in equation \eqref{derivative term} for the independent space-time random walk.
\qed
\end{proofsect}

\begin{proofsect}{Proof of Proposition~\ref{diagonalization_theorem}}
We only sketch the proof of the Markovian case, the Bosonic one follows analogously. The independence follows from Campbell's formula: Since the Laplace transform of the sum of the two Poisson point processes is the exponential of the sum of two integrals with respect to the intensity measures, the result can be factorised. To calculate the distribution of the point loops one makes the observation
\begin{equation}\label{pointloops}
\sum_{x\in\gb}\int_0^\infty\left( 1-\ex^{-\langle v,L\rangle}\right)\frac{\ex^{\mu t}}{t}\P_{x,x}^{\ssup{t}}(\{\text{ no jump up to time }t\})\d t=\sum_{x\in\gb}\int_0^\infty\left(1-\ex^{-v(x)t}\right)\frac{\ex^{-t(-d(x)-\mu)}}{t}\d t\, .
\end{equation}
From then on, one can proceed as in the proof of Proposition \ref{properties}.\qed
\end{proofsect}

\begin{proofsect}{Proof of Example~\ref{counterexample}}
We show that  our convergence results in the torus limit $ N\to\infty $  fail to hold (Theorem~\ref{THM-main1} and Theorem~\ref{THM-main2}) when the space-time random walk has equal rates for going up or down along the ``time'' torus $ \T_N$. An easy computation shows that symmetric rates along the torus imply that in the equation after \eqref{loop_measure_convergence_k=1} the factor $\P^{\ssup{t}}_{N,0,0}$ is replaced by its square  $\big(\P^{\ssup{t}}_{N,0,0}\big)^2$, as we can circle around the torus in two different directions. It  follows from the previous proofs  that for all $t>0$ the convergence of $\P^{\ssup{t}}_{N,0,0}\to 0$ is uniform  in $t$ at least of order square-root in $N$. Actually, for finite $ t$ the term vanishes exponentially in $N$.  Therefore all finite dimensional distributions vanish in the torus limit and thus the limiting procedure does  not establish a loop measure for the projected space-time process.\qed
\end{proofsect}

\bigskip

\begin{proofsect}{Proof of Theorem~\ref{support_prop}}
The spatial projection of the given Bosonic walk is denoted $M_N $. We abbreviate $f(t)=\1_A(t)$. In the definition of the projected loop measure $M_N $, we replace $ \ex^{\mu t}/t $ by $ f(t) $ (and we can later employ an approximation scheme) as that factor does not change the support property  of the loop measure in the torus limit. We only need to consider paths which wind around the torus once as all other paths can be written as a concatenation of such paths. Thus we shall estimate the probability that such loops have length in the given set $A$. We will then use that result to prove the theorem.
Denote $ \widetilde{d}(x,\tau) $ the diagonal elements of the the generator $G_N $.  We calculate the weight of all space-time paths from $(x,\tau) $ to $ (y,\tau) $ winding around the torus only once:
\begin{equation}\label{long support prop expand proof}
\begin{split}
&\sum_{\tau=0}^{N-1}\int_0^\infty f(t)\P_{(x,\tau)}(X_t=(y,\tau), \wind(X_{[0,t]})=1)\d t\\
&=\frac{1}{N}\sum_{\tau=0}^{N-1}\sum_{\heap{z\colon (x,\tau)\to(y,\tau),}{\wind(z)=1}}\int_0^\infty f(t)\int_0^t \ex^{(t-s)\widetilde{d}(y,\tau)}N P(z)\d \P(W(z)=s)\d s\, \d t\\
&=\frac{1}{N}\sum_{\tau=0}^{N-1}\sum_{\heap{z\colon (x,\tau)\to(y,\tau),}{ \wind(z)=1}}P(z)\E_z\Big[N\int_X^\infty f(t)\ex^{(t-X)\widetilde{d}(y,\tau)}\d t\Big]\\
&=\frac{1}{N}\sum_{\tau=0}^{N-1}\sum_{\heap{z\colon (x,\tau)\to(y,\tau),}{ \wind(z)=1}}P(z)\int_0^\infty N\ex^{t\widetilde{d}(y,\tau)}\E_z\Big[f(t+X)\Big]\d t\\
&=\frac{1}{N}\sum_{\tau=0}^{N-1}\sum_{\heap{z\colon (x,\tau)\to(y,\tau),}{ \wind(z)=1}}P(z)\Ex_{-\widetilde{d}(y,\tau)}\Big[\frac{N}{-\widetilde{d}(y,\tau)}\E_z\Big[f(\texttt{C}+X)\Big]\Big]
\, ,
\end{split}
\end{equation}
where $ \E_z $ denotes the expectation for the random waiting time $W(z)$ of the given path $=z$. $ \Ex_{-\widetilde{d}(x,\tau)} $ denotes the expectation with respect to an exponential distribution with parameter $ \widetilde{d}(x,\tau) $ for the random time $ \texttt{C} $. 
Firstly, assume that $\beta\in (\R\setminus A)^\circ$. Then, for some $\epsilon >0$ it holds that
\begin{equation}\label{L2 bound}
\P_z(W\in A)\le \P_z(|W-\beta|>\epsilon)\le\epsilon^{-2}\E_z\left[\left(W-\beta\right)^2\right]\, .
\end{equation}
Expanding the square and calculating the expectations, and writing 
$$ 
\d z = \frac{1}{N}  \sum_{\tau=0}^{N-1}\sum_{\heap{z\colon (x,\tau)\to(y,\tau),}{ \wind(z)=1}}P(z),
$$ using the assumptions from Theorem \ref{support_prop} and equations \eqref{long support prop expand proof} together with \eqref{L2 bound}, one sees that
\begin{equation}
\lim_{N\to\infty}\int\d z\P_z(W\in A)=0\, ,
\end{equation}
which was the claim. If $\beta\in A^\circ$, we simply use
\begin{equation}
\int\d z\P_z(W\in A)=\int\d z -\int\d z \P_z(W\in\R\setminus A)\, ,
\end{equation}
and the apply the same reasoning as above to obtain
\begin{equation}
\lim_{N\to\infty}\int\d z\P_z(W\in A)=d_{x,y}>0\, .
\end{equation}
One can now approximate $\ex^{\mu t}/t$ by step functions to conclude the argument.
\medskip

\noindent  We now prove the last statement, namely that all weakly and strongly asymptotically independent space-time random walks are Bosonic random walks. Note, all strongly asymptotically independent space-time random walks are also weakly asymptotically independent. The difficulty to establish that weakly asymptotically independent random walks belong to the Bosonic loop walk class comes from the fact that the jump rates on $\Lambda$ are not uniform.  
We first show \eqref{waitprob} and \eqref{waitingtime} for the independent space-time random walk with uniform jump rates. To establish \eqref{waitprob} we can use the independence to obtain from equation \eqref{step probabolities}, 
\begin{equation}
P(z)=\Big(1-\frac{\lambda\beta}{N+\beta\lambda}\Big)^N\Big(\frac{\lambda}{\lambda+N\beta^{-1}}\Big)^kp(x,x_1)\cdots p(x_{k-1},y)\, .
\end{equation}
Thus, we write
\begin{equation}
\frac{1}{N}\sum_{\tau=0}^{N-1}\sum_{\heap{z\colon (x,\tau)\to(y,\tau)}{\wind(z)=1}}P(z)=\Big(1-\frac{\lambda\beta}{N+\beta\lambda}\Big)^N\sum_{k=0}^\infty\binom{N+k}{k}\Big(\Big[\frac{\lambda}{\lambda+N\beta^{-1}}\Big]P\Big)^k(x,y)\, ,
\end{equation}
which obviously has a well defined limit
$$
\lim_{N\to\infty} \frac{1}{N}\sum_{\tau=0}^{N-1}\sum_{\heap{z\colon (x,\tau)\to(y,\tau)}{ \wind(z)=1}}P(z)\in (0,1].
$$ 
To check the second condition \eqref{waitingtime}, we rewrite it as
\begin{equation}\label{eq429}
\sum_{k=0}^\infty\binom{N+k}{k}\Big(\Big[\frac{\lambda}{\lambda+N\beta^{-1}}\Big]P\Big)^k(x,y)\Big[\Big(\beta\Big(1-\frac{N}{N+\beta}\Big)-\frac{k}{N\beta^{-1}+\lambda}\Big)^2	+\frac{N+k}{(N\beta^{-1}+\lambda)^2}\Big]\, ,
\end{equation}
which converges to zero since $P^k(x,y)$ decays exponentially in $k$. If the jump rates are not uniform consider a path as above. Then
\begin{equation}
\begin{split}
\sum_zP(z)=&p(x,x_1)\frac{\lambda_x}{N\beta^{-1}+\lambda_x}\cdots p(x_{k-1},y)\frac{\lambda_{x_{k-1}}}{N\beta^{-1}+\lambda_{x_{k-1}}}\sum_{\heap{j_0,\ldots j_k}{ j_0+\ldots+j_k=N}}\prod_{i=0}^k\Big(\frac{N\beta^{-1}}{N\beta^{-1}+\lambda_{x_i}}\Big)^{j_i}\, ,
\end{split}
\end{equation}
where the sum is over all paths on the space-time graph which have jumps $x=x_0\to x_1\to\ldots x_{k-1}\to x_k=y$ on the lattice. As we now sum over all the different $x_1,\ldots,x_k$,  we get
\begin{equation}
\begin{split}
&\sum_{z,x_1,\ldots,x_k}P(z)\\
&=\sum_{x_1,\ldots,x_k}p(x,x_1)\frac{\lambda_x}{N\beta^{-1}+\lambda_x}\cdots p(x_{k-1},y)\frac{\lambda_{x_{k-1}}}{N\beta^{-1}+\lambda_{x_{k-1}}}\sum_{\heap{j_0,\ldots j_k}{ j_0+\cdots+j_k=N}}\prod_{i=0}^k\Big(\frac{N\beta^{-1}}{N\beta^{-1}+\lambda_{x_i}}\Big)^{j_i}\\
&=\sum_{x_1,\ldots,x_k}p(x,x_1)\frac{\lambda_x}{N\beta^{-1}+\lambda_x}\cdots p(x_{k-1},y)\frac{\lambda_{x_{k-1}}}{N\beta^{-1}+\lambda_{x_{k-1}}}\sum_{i=0}^k\frac{\Big(1-\frac{\beta \lambda_{x_i}}{N+\lambda_{x_i}\beta}\Big)^{k+N-1}}{\prod_{\heap{l=0}{ l\neq i}}^k\Big[\frac{\beta \lambda_{x_l}}{N+\lambda_{x_l}\beta}-\frac{\beta \lambda_{x_i}}{N+\lambda_{x_i}\beta}\Big]}\\
&\to\sum_{x_1,\ldots,x_k}p(x,x_1)\lambda_x\beta\cdots p(x_{k-1},y)\beta \lambda_{x_{k-1}}\sum_{i=0}^k\frac{\ex^{-\beta \lambda_{x_i}}}{\prod_{\heap{l=0}{ l\neq i}}^k\big[\beta \lambda_{x_l}-\beta \lambda_{x_i}\big]}\, \; \mbox{ as } N\to\infty\,,
\end{split}
\end{equation}
where we have used Theorem 3.2 from \cite{C11} (or simply using partial fraction decomposition). If not all $\lambda_{x_i}$'s are distinct,  the above formulae hold in their canonical extension. Similarly, we show that the second condition is satisfied, i.e., the analogue of equation \eqref{eq429} converges to zero. To conclude the proof, we show that every path generated by weights which are zero for the independent case but  that do not vanish for the perturbed case have exponentially small contributions.  
Fix some $\epsilon>0$. In the expression \eqref{waitprob}  paths can, a priori, have a length of $N+k$ with $k\ge -N$. As we assume that  $\kappa(x)-\mu>0$ for at least one $x\in\Lambda$, we can choose some $M\in\N$ such that paths with length greater than $N+M$ have weight less than $\epsilon$. We thus restrict ourselves to paths of length smaller than $N+M$.\\
Given a path of length $A+B\le N+M$, we can decompose it into $A$ jumps along edges $(z_i,z_{i+1})$ for which $Q_N(z_i,z_{i+1})>0$ and $B$ jumps for which the $Q_N$ entries are zero. Fixing the $A$ part of the path, there are less than $(|\Lambda|N)^{AB}\le(|\Lambda|N)^{NB+MB}$ paths with $B$ jumps. Each of these paths has a weight dominated by $N^{-2N|\Lambda|B}$. Summing over all $B$, we can see that the contribution from all perturbed paths vanishes. Similarly, given a path for which $B=0$, we write the weight of each jump as $Q_N+\Ocal(N^{-2|\Lambda|N})$. Expanding the product, we again see that the contributions are negligible. As $\epsilon$ was arbitrary, the result follows.
\qed
\end{proofsect}

\bigskip

\begin{proofsect}{Proof of Example~\ref{ergodic example}}
The probability density of the space-time random walk with periodic boundary conditions is 
\begin{equation}
p^{\texttt{per}}_t(x,0,y,\tau)=\ex^{-t\beta^{-1}N}\sum_{k=0}^\infty\frac{(t\beta^{-1}N)^{Nk+\tau}}{(\tau+Nk)!}P^{Nk+\tau}(x,y)=\E[P^{S_N+\tau}(x,y)\mathbbm{1}\{S_N=\tau\!\!\! \mod N\}]\, ,
\end{equation}
$ x,y\in\L_M, \tau\in\T_N $,
where $S_N$ is Poisson distributed with parameter $\beta^{-1}N$. Furthermore,
$$
\begin{aligned}
&\E[P^{S_N+\tau}(x,y)\1\{S_N=\tau\!\!\! \mod N\}]=\pi(y)\E[\1\{S_N=\tau\!\!\! \mod N\}]\\[1.5ex]& + \E[(P^{S_N+\tau}(x,y)-\pi(y))\1\{S_N=\tau\!\!\! \mod N\}]\,.
\end{aligned}
$$
The assumptions for the spatial component ensure that there is a $c>0$ such that for all $x,y$ and $N$ large enough, 
\begin{equation}
|P^N(x,y)-\pi(y)|\le \ex^{-Nc}\,.
\end{equation} 
Therefore, 
\begin{equation}
p^{\texttt{per}}_t(x,0,y,\tau)=\pi(y)\P^{\ssup{t}}_{N,0,\tau}+\Ocal(\ex^{-cN})\,
\end{equation}
as
$$
\E[(P^{S_N+\tau}(x,y)-\pi(y))\1\{S_N=\tau\!\!\! \mod N\}]  \le \E[ |P^{S_N}(x,y)-\pi(y)|]\le \E[\ex^{-cS_N}]\le \ex^{-cN}\, ,
$$
From then on, the result follows repeating the steps  in the proof of Theorem \ref{THM-main1}.
\qed
\end{proofsect}

\subsection{Isomorphism theorems}\label{Proof2}

The key identity to the proof of Theorem \ref{Isomorphism occupation field} is stated as a separate lemma due to its importance.
\begin{lemma}\label{laplacetransform}
Assume that  $\kappa(x)-\mu>0$ for at least one $x\in\Lambda$. Let $v\in [0,\infty)^\gb$. Then 
\begin{equation}\label{laplacetransformeq}
\Er_{\gb,\mu}[\ex^{-\langle v,\Lscr \rangle}]=\frac{\det(Q+\mu I)}{\det(Q+\mu I-V)}\, ,
\end{equation}
where $V=\diag(v)$ denotes the diagonal matrix with $v$ on the diagonal.
\end{lemma}

\begin{proofsect}{Proof of Lemma \ref{laplacetransform}}
We use Campbell's formula to rewrite the left-hand side of \eqref{laplacetransformeq} as
\begin{equation}
\Er_{\gb,\mu}\Big[\ex^{-\langle v,\Lscr\rangle}\Big]=\exp\Big\{-\int_\Gamma\Big(1-\ex^{-\langle v,L\rangle}\Big)\d M_{\gb,\mu}\Big\}\, .
\end{equation}
With the help of the Feynman-Kac formula one can now write
\begin{equation}
\begin{aligned}
\int_\Gamma\Big(1-\ex^{-\langle v,L\rangle}\Big)\d M_{\gb,\mu}&=\sum_{x\in\gb}\int_0^\infty\frac{\ex^{t\mu}}{t}\Big(\ex^{tQ}(x,x)-\ex^{t(Q-V)}(x,x)\Big)\d t\\[1.5ex]& =\int_0^\infty\frac{\ex^{t\mu}}{t}\Big(\text{Tr}[\ex^{tQ}]-\text{Tr}[\ex^{t(Q-V)}]\Big)\d t.
\end{aligned}
\end{equation}
One then uses
\begin{equation}
\ex^{t\mu}\text{Tr}[\ex^{tQ}]=\sum_{x\in\gb}\ex^{t(q(x)+\mu)}\, ,
\end{equation}
where $q(x)$ is the $x$'th eigenvector of $Q$, analogously $q^v(x)$ for $Q-V$. Hence, using again a version \eqref{logformula}, 
\begin{equation}
-\int_\Gamma\Big(1-\ex^{-\langle v,L\rangle}\Big)\d M_{\gb,\mu}=\sum_{x\in\gb}\log\Big(\frac{q(x)+\mu}{q^v(x)+\mu}\Big)=\log\Big(\frac{\det(Q+\mu I)}{\det(Q+\mu I-V)}\Big)\, .
\end{equation}
\qed
\end{proofsect}

The Laplace transform of the Bosonic Loop measure is computed as follows.
\begin{lemma}
For $v\in[0,\infty)^\gb$ we have that
\begin{equation}
\Er_{\gb,\mu,\beta}^B\Big[\ex^{-\langle v,\Lscr\rangle}\Big]=\frac{\det\Big(I-\ex^{\beta(Q+\mu I)}\Big)}{\det\Big(I-\ex^{\beta(Q+\mu I-V)}\Big)}\, .
\end{equation}
\end{lemma}
\begin{proofsect}{Proof} 
One expands the expectation analogously to the Markovian case and uses  the Campbells formula again to calculate the expectation. Using the Feynman-Kac formula one  calculates the expectation of the local time and then uses the reasoning from the proof of Proposition \ref{properties}.
\qed
\end{proofsect}

\bigskip

\noindent We now turn to the proof of the Isomorphism theorems.
\begin{proofsect}{Proof of Theorem~\ref{Isomorphism occupation field}}
(a) Our assumptions ensure that $C=A^{-1}  $ and thus $A$ has a positive Hermitian part, henceforth $ \mu_A $ is the complex Gaussian measure on $ \C^\gb $. 
We shall compute Laplace transforms separately on both sides of equation \eqref{Isoa} , that is, for the measure $ \mu_A $ and for the Poisson measure $ \Pr_{\gb,\mu} $. For that, let 
$$ 
v_0\in \R^{\gb} :=\inf_{v\in\R^{\gb}}\big\{\E_A\big[\ex^{-\langle\phi,V\overline{\phi}\rangle}\big]<\infty\big\}\,,
$$
where $ V =\diag(v) $ denotes the diagonal matrix with entries $v$ and where we order the fields component-wise. Gaussian calculus thus ensures that the Laplace transform for the measure $ \mu_A$ is given
$$
\E_A[ \ex^{-\langle \overline\phi,V\phi\rangle}]=\frac{Z_{A+V}}{Z_{A}}=\frac{\det(Q+\mu I)}{\det(Q+\mu I-V)}\; \mbox{ for all }\;   V\ge V_0.
$$
In fact, we can set $ v_0 \equiv 0 $ here as the occupation field assumes only positive values.
The Laplace-transform for the Poisson measure $ \Pr_{\gb,\mu} $ has been calculated  in Lemma \ref{laplacetransform},
\begin{equation}\label{laplacetransform_equation}
\Er_{\gb,\mu}[\ex^{-\langle v,\Lscr\rangle}]=\frac{\det(Q+\mu I)}{\det( Q+\mu I-V)}\, .
\end{equation}
Hence, the Laplace transforms for both measures agree.  The family of functions $t\mapsto \ex^{-\langle v,t\rangle }$ for $v\in [0,\infty)^\L$ is a point separating family. Having shown that for such functions, the expectation of the square of the field  equals the expectation of the occupation field, the  two measures agree. As the Poisson point measure is a probability measure, the square of the complex Gaussian field is governed by a probability measure as well.

\medskip

\noindent (b)  We now prove the second identity of Theorem \ref{Isomorphism occupation field}. Define the complex ``density" $f_\C$ (Radon-Nikodym derivative) by
\begin{equation}
f_\C(\phi)=\frac{1}{Z_S}\ex^{-\langle \phi,A^S\bar\phi\rangle}\, ,
\end{equation}
where $A^S$ is the skew Hermitian part of $A$ and $ Z_S=\frac{\det(\widetilde{A})}{\det(A)} $.
Then, for any bounded measurable function $F\colon \C^\gb\to\C$, 
$$
\E_A[F(\phi)]=\E_{\widetilde{A}}\left[f_\C(\phi)F(\phi)\right]\, .
$$
Note that as $\widetilde{A} $ is Hermitian, the resulting measure $\mu_{\widetilde{A}} $ on $\C^\gb$ is a probability measure. This allows, for any sub-$\sigma$ algebra $\Fcal$ on $\C^\gb$ to define the ``conditional expectation"
\begin{equation}
\E_A[F(\phi)|\Fcal](\tilde\phi)=\E_{\widetilde{A}}[f_\C(\phi)F(\phi)|\Fcal](\tilde\phi)\, .
\end{equation}
Note that this implies the tower property for the conditional expectation.\\
Using that conditional expectation is integration with respect to the regular conditional probability, we can apply this to the sigma algebra $\Fcal_|$ generated by $\phi\mapsto|\phi|$ to get
$$
\E_A[F(\phi)|\Fcal_|](\tilde\phi)=\int f_\C(\phi)F(\phi)\mu_{\widetilde{A}}[\phi\in\d\phi|\Fcal_|](\tilde\phi)=\int f_\C(\phi)F(\phi)\mu_{\widetilde{A}}[\phi\in\d\phi||\phi|=\tilde\phi]\, , 
$$
$ \tilde\phi\in[0,\infty)^\gb\, $, where $\mu_{\widetilde{A}}[\phi\in\d\phi||\phi|=\tilde\phi]$ denotes the regular conditional probability kernel associated to $\Fcal_|$.
We can compute $\mu_{\widetilde{A}}[\phi\in\d\phi||\phi|=\tilde\phi]$, since,  by parametrising $\C^\gb$ using polar coordinates,  the joint density of $(\theta,|\phi|)$ is given by
\begin{equation}
\mu_{\widetilde{A}}\big(\{(\theta,|\phi|)\in \Acal\}\big)=\frac{1}{Z_H}\int_{\mathbb{S}^\gb\times[0,\infty)^\gb}\mathbbm{1}_{\Acal}(\theta,|\phi|)\ex^{-\langle |\phi|\theta,\widetilde{A}|\phi|\bar{\theta}\rangle}\underbrace{\Big[\prod_{x\in\gb}\d \omega_x(\theta)\Big]}_{:=\d\Scal_\gb}\Big[\prod_{x\in\gb}\d|\phi(x)|\Big]\, ,
\end{equation}
for some event $ \Acal $,  where $\d \abs{\phi(x)}$ is integration with respect to  the Lebesgue measure on $[0,\infty)$ and $\omega_x$ is the surface measure on the unit sphere $\mathbb{S}=\{z\in\C\colon |z|=1\}$.
By Bayes Theorem, the density of the conditional probability is
$$
\mu_{\widetilde{A}}[\phi\in\d\phi||\phi|=\tilde\phi]=\frac{\ex^{-\langle \tilde\phi\theta,\widetilde{A} \tilde\phi \bar{\theta}\rangle}}{\int\ex^{-\langle \tilde\phi\theta,\widetilde{A} \tilde\phi \bar{\theta}\rangle}\d \Scal_\gb(\theta)}\, ,
$$
where we have implicitly identified $\phi=\theta \tilde\phi $.  We now apply (a) to the square $ \widetilde{\phi}^2=(\sqrt{\tilde\phi}_x\sqrt{\tilde\phi}_x)_{x\in\gb} $, and obtain
$$
\begin{aligned}
\E_A[F(\phi)]&=\E_A\big[\E_A\big[F(\phi)\big|\Fcal_|\big]\big]=\E_{\widetilde{A}}\big[\E_{\widetilde{A}}\big[f_\C(\phi)F(\phi)|\Fcal_|\big]\big]
\\[1.5ex]&=\Er_{\widetilde{A}}\big[\int\, f_\C(\theta\sqrt{\Lscr})F(\theta\sqrt{\Lscr})\ex^{-\langle \sqrt{\Lscr}\theta,\widetilde{A}\sqrt{\Lscr}\bar{\theta}\rangle}\d\Scal_\gb(\theta)\big]\, .
\\[1.5ex]& = \Er_{\widetilde{A}}\big[\frac{1}{Z_S\int\ex^{-\langle \sqrt{\Lscr}\theta,\widetilde{A}\sqrt{\Lscr}\bar{\theta}\rangle}\d \Scal_\gb(\theta)}\int F(\theta\sqrt{\Lscr})\ex^{-\langle \sqrt{\Lscr}\theta,A\sqrt{\Lscr}\bar{\theta}\rangle}\d\Scal_\gb(\theta)\big]\,.
\end{aligned}
$$

\medskip

\noindent (c) To prove the identity for non-symmetric matrices $A$, we follow \cite{BIS09}. Let us split $A$  into its symmetric (or Hermitian) $\widetilde{A}$ and skew-symmetric (skew Hermitian) part $A^S$, namely define 
\begin{equation}\label{matrixID}
A(z)=\mu \1+\widetilde{A}+zA^S\,,z\in\C.
\end{equation}  Our goal is to show that the identity \eqref{matrixID} holds for $z=1$ as this gives us back $A$. From the previous computation we know \eqref{matrixID}  holds for $z=0$. We shall show the identity \eqref{matrixID} for purely imaginary $z=ri$ with  $0<r<\epsilon$ for some $\epsilon>0$ sufficiently small. However, for such $z$'s the generator would be Hermitian but also complex valued. We thus need to examine whether  the previous results carry over to the complex case. We  use the same notation as in the non-complex case.\\ One  can define the complex valued measure $\P_x$ and its associated semigroup $S(t)$ through its transition density
\begin{equation}
p_t(x,y)=\ex^{-tA(z)}(x,y)\, ,
\end{equation}
and expanding this to all cylinder sets on the Skorohod space $ D_\gb $ and using this density to obtain  a measure on the full space $D_\gb$. Note that this is no longer a Markov-semigroup (as positivity preservation is lost) but for $\epsilon$  sufficiently small  it still establishes a uniformly continuous $C_0$ semigroup. In particular $\partial_tS(t)=-A(z)S(t)=-S(t)A(z)$. \\
To show that the distribution of $(\bar{\phi_x}\phi_x)_{x\in\gb}$ agrees with the one of $(\Lcal_x)_{x\in\gb}$ we only used two probabilistic tools:  Campbell's formula and the Feynman-Kac formula. As the proof of Campbell's formula does not use the fact that the intensity measure is positive, one can transfer it to the complex setting.\\
For $v\in\C^\gb$ with $ \norm{v} $ ($V=\diag(v) $ denotes the diagonal matrix with entries $ v$) sufficiently small,  the Feynman-Kac formula allows to define  the Feynman-Kac semi-group $\texttt{FK}_V(t)$  as
\begin{equation}
\texttt{FK}_V(t)f(x)=\E_x\Big[f(X_t)\ex^{-\int_0^tV(X_S)\d s}\Big]\, ,
\end{equation}
has derivative $\partial_t\texttt{FK}_V(t)=(-A(z)+V)\texttt{FK}_V(t)$ and thus 
\begin{equation}
\texttt{FK}_V(t)f(x)=\ex^{t(-A(z)+V)}f(x)\, .
\end{equation}
The finiteness of $M_{\Lambda,\mu}\Big[1-\ex^{-\langle v,L\rangle}\Big]$ is still given as $A(z)$ remains to have positive eigenvalues. The span of the family
\begin{equation}
\{\ex^{-\langle v,|\phi|^2\rangle}\colon v_x\; \ge 0\, \forall x\in\gb\}\, ,
\end{equation}
is a point separating family (over $[0,\infty)^\gb$),  and by the Stone-agreedWeierstrass Theorem it suffices to compare measures on this class. By the previous result  for $\epsilon>0$ sufficiently small
\begin{equation}
\E_{A(z)}[\ex^{-\langle v,\abs{\phi}^2\rangle}]=\Er_{A(z),\gb}[\ex^{-\langle v,{\Lcal}\rangle}]=\frac{\det(A(z)-\mu I)}{\det(A(z)-\mu I+V)}\, .
\end{equation}
Using that both sides are meromorphic in $z$ (this is due to the definition of the determinants, see \cite{BIS09}), one can extend the result to the whole complex plane and in particular to the desired case $z=1$.\\
Note that for $z$ imaginary and small, both sides are probability measures, despite the matrix $A(z)$ being complex, this is due to the fact that $A(z) $ is  Hermitian. For probability measures one  can use the standard theory of conditional expectations to verify the isomorphism equation in (b) above. Then one chooses a point separating family over $\C^\gb$ to show analyticity in $z$. Here, the polynomials $P(\phi)=\left(\prod_i\phi(x_i)\right)\left(\prod_j\overline{\phi(y_i)}\right)$ are suitable,  and their analyticity follows from the calculations performed in \cite{BIS09}. Via the given  density we extend the result to all $F$ measurable and bounded, i.e.,
\begin{equation}
\E_A[F(\phi)]=\Er_{\gb,\mu}\Big[\frac{1}{\int \ex^{-\langle\theta\sqrt{\Lscr},A\sqrt{\Lscr}\bar{\theta}\rangle}\d \Scal_\gb(\theta)}\int F(\theta\sqrt{\Lscr})\ex^{-\langle\theta\sqrt{\Lscr},A\sqrt{\Lscr}\bar{\theta}\rangle}\d \Scal_\gb(\theta)\Big]\, .
\end{equation}
As a by-product of this proof, we get that the previous identities hold for complex weighted loops in the above fashion.

\qed

\end{proofsect}

\begin{proofsect}{{Proof of Proposition~\ref{Symanzik_formula_generalized}}}
We set $\beta=1$, as the other values of $ \beta $ can be handled by change of variables, and let $ \mu\le 0 $.  Note that by combining Dynkin's Isomorphism theorem and Le Jan's Isomorphism (Theorem \ref{Isomorphism occupation field}) we have the new isomorphism theorem
\begin{equation}\label{newisomorohism}
\E_{xy}^{\mu,1}\otimes\Er_{\gb,\mu}[J(\Lscr+L)]=\E_A[\overline{\phi}_x\phi_yJ(|\phi|^2)]\, , \mbox{ for any bounded and continuous } J\colon[0,\infty)^\gb\to\R.
 \end{equation}
We can now equate
\begin{equation}
\begin{split}
\E_{A,J}[\overline{\phi}_x\phi_y]&=\frac{1}{Z_{A,J}}\left[Z_A\E_A\left[\overline{\phi}_x\phi_yJ(|\phi|^2)\right]\right]\\[1.5ex]
&=\frac{Z_A}{Z_{A,J}}\E_{xy}^{\mu,1}\otimes\Er_{\gb,\mu}[J(\Lscr+L)]\, .
\end{split}
\end{equation}
The denominator is computed analogously. For higher moments, one notes that Dynkin's Isomorphism theorem holds not only for two point correlations but can be extended to arbitrary $k$-point correlations.
\qed
\end{proofsect}

\bigskip

\begin{proofsect}{Proof of Theorem \ref{Space-time loop-representation of the correlation function}}
(a) (i) This can actually be done similarly to the proof of Theorem~\ref{THM-main2}, that is, first by proving the statements for the independent space-time random walk with generator $G_N=Q_N $ (see Definition~\ref{definition of the walk}), and in a second step to extend all statements to generators $G_N $ of weakly asymptotically independent space-time random walks.  However, there is one  difference: we cease to have the factor $\tfrac{1}{t}$ in the definition of the bridge measure which causes a divergence at zero. We obtain as in \eqref{splitting},
$$
\begin{aligned}
\E_{x,y,B}^{\mu,\beta}[F(L)]=\sum_{x\in\L}\sum_{j=0}^\infty\int_0^\infty\,\ex^{\beta\mu t} N\E_{x,y}^{\ssup{\beta t}}\big[F((L_x)_{x\in\L})\big]\P_{N,0}^{\ssup{\beta t}}(X_t^{\ssup{2}}=0,\wind =j)\,\d t,
\end{aligned}
$$ 
and it is straightforward to see that the part with  winding numbers $ \wind\ge 1 $ does converge to the Bosonic path measure (in fact due to the missing singularity the proof here is easier than the one in Theorem~\ref{THM-main2}). The term with $ \wind=0 $ follows as
$$
\begin{aligned}
\int_0^\infty & \,\ex^{\beta \mu t} N\ex^{-t N}\Big(F(t_x)\ex^{td_x}\delta(x,y)+(1-\ex^{td_x})\E_{x,y}^{\ssup{\beta t}}\big[F(L)\1\{\mbox{at least one jump}\}\big]\Big)\\[1.5ex]
&=\tilde\E_{N}\Big[ \ex^{\beta X} \Big(F(X)\ex^{Xd_x}\delta(x,y)+(1-\ex^{ Xd_x})\E_{x,y}^{\ssup{\beta X}}\big[F(L)\1\{\mbox{at least one jump}\}\big]\Big)\Big]\\[1.5ex]
& \; \longrightarrow  F(0)\delta(x,y) \; \mbox{ as } \; N\to\infty,
\end{aligned}
$$
where $ \tilde\E_N $ is expectation with respect to an exponentially distributed random variable with expectation $1/N $.

\noindent (ii) Clearly, Theorem~\ref{THM-main2} implies that
$$
\lim_{N\to\infty}\Er_N[J(\Lscr)]=\Er^B_{\L,\mu,\beta}[J(\Lscr+\beta)].
$$
This implies \eqref{RDMresult1} using Theorem~\ref{physical quantities representation}. To see that this stems from Symanzik's formula (Proposition~\ref{Symanzik_formula_generalized}) in the space-time setting, we sum equation \eqref{momentST} over all possible torus points and note that $ J(L+\Lscr)=\Jsf(\psf_\L(\Lb+\bLscr)) $, i.e.,
$$
\begin{aligned}
\sum_{\tau=0}^{N-1}\E_{A_N,\Jsf}\big[\overline{\bphi}_{(x,\tau)}\bphi_{(y,\tau)}\big]&=\sum_{\tau=0}^{N-1} \frac{1}{\Er_{\L\times\T_N,\mu}[\Jsf(\bLscr)]}\E_{(x,\tau),(y,\tau) }^{\mu,\beta}\otimes\Er_{\L\times\T_N,\mu}\big[\Jsf(\Lb+\bLscr)\big]\\[1.5ex]
&=\frac{\E_{x,y,N}^{\mu,\beta}\otimes \Er_N[F(L+\Lscr)]}{\Er_N[F(\Lscr)]}.
\end{aligned}
$$

\noindent (b) This follows easily from (a). The constant background field $ \beta $ drops out due to the linearity of the interaction. 
\qed
\end{proofsect}

\bigskip

\begin{proofsect}{Proof of Theorem \ref{infvollimprop}}
(a) The statements in Proposition~\ref{properties} carry over as the finiteness of the graph is not used in the proof.  (b) and (c) follow by taking advantage of the convergence in the transient regime for $ d\ge 3 $ and $ \mu\le 0 $ or $ d\ge 1 $ and $ \mu<0 $. For (c) we require the support of the functional $ F $ to be compact to ensure finiteness when passing to the thermodynamic limit. (d)  In Theorem 2.7 we only consider point loops. For any finite $v\in [0,\infty)^{\Z^d}$ with $v(x)=0$ outside some finite set $\Lambda$, their distribution is identical with the distribution of the point loops under the loop measure $ M_{\L,\mu} $ or the Bosonic loop measure $M_{\Lambda,\mu,\beta}^B$, as these point loops do not leave the set $\Lambda$. As the Poisson process of loops forms a consistent family of measure in the sense of Kolmogorov's extension theorem (see \cite{Sznitman}), the class of such  $v$ is sufficient. With \eqref{pointloops} we conclude with the statement.  \\
\noindent (e) We sketch the proof for extension to $ \Z^d $ of the statement in (a) of Theorem~\ref{Isomorphism occupation field}. Statement (b) is included in (a) as the generators are symmetric, and (c) can be obtained using characteristic functions instead of Laplace transforms. We are using Laplace transforms for all distributions of $ |\phi|^2=(\overline{\phi}_x\phi_x)_{x\in\L} $ in (a).   
We denote  $c_0(\Z^d)$ the set of all vectors $v$ in $\R^{\Z^d}$ such that all but finitely many coordinates are zero. 
Define the set 
\begin{equation}
\{f_v(x)=\ex^{-\langle v,x\rangle}\colon v\in c_0(\Z^d)\text{ and } 0\le v(x)<\infty\text{ for all }x\in\Z^d\} \, ,
\end{equation}
of functions on $ \Z^d $. This set is a point-separating family stable under multiplication and including the constant function 
for  continuous and bounded functions from $[0,\infty)^{\Z^d}$ to $\R$. We apply the one point compactification to enlarge the space to $[0,\infty]^{\Z^d}$ which is compact by Tychonoff's theorem. We set for $g\colon [0,\infty)\to\R$, 
\begin{equation}
g(\infty)=\lim_{t\to\infty}g(t)\, ,
\end{equation}
and continue this component wise to functions from $[0,\infty)^{\Z^d}$, requiring that the limit is well defined. By the Stone-Weierstrass theorem, it suffices to check that all Laplace transforms, i.e. expectations of  $(f_v) $ for any $ v\in c_0(\Z^d) $,  of two measures agree to conclude that they coincide. We compute  (recall that $ A_\L=-Q_\L $ with $Q_\L-\mu\1_\L$ being the generator of the random walk on $ \L $ and that $ V $ is the diagonal  matrix with entries $v(x) $)
\begin{equation}
\E_{A_\Lambda}\big[\ex^{-\langle \phi,V\overline{\phi}\rangle}\big]=\frac{\det(A_\L)}{\det(A_\L+V)}= \det(\1_\L+A_\L^{-1}V)^{-1}=\det(\1_\L+G_\L^\mu V)^{-1},
\end{equation} 
where the Green function $ $ is defined as 
\begin{equation}
G_\L^\mu(x,y)=\int_0^\infty\,\ex^{\mu t}\,p_t(x,y)\,\d t\,,\quad x,y\in\L.
\end{equation}
For $ d\ge 1 $ and $ \mu<0 $, or for $d\ge 3 $ and $ \mu\le 0 $, the Green function is finite (transient). Therefore  the limit of the determinants exists  as a Fredholm determinant, and using \eqref{laplacetransform} implies
\begin{equation}
\lim_{\Lambda\uparrow\Z^d}\E_{A_\Lambda}\big[\ex^{-\langle \phi,V\overline{\phi}\rangle}\big]=\det(\1_{\Z^d}+G^\mu_{\Z^d}V)^{-1}=\Er_{\Z^d,\mu}\big[\ex^{-\langle v,\Lscr\rangle}\big]=\lim_{\L\uparrow\Z^d}\Er_{\L,\mu}\big[\ex^{-\langle v,\Lscr\rangle}\big]\, ,
\end{equation}
and thus both distributions agree. We are using  that the distribution $\Pr_{\Lambda,\mu}\circ\Lscr^{-1} $ converges to the distribution $\Pr_{\Z^d,\mu}\circ\Lscr^{-1} $ with respect to  the topology of local convergence, which can be obtained from the convergence of the  intensity measure $M_{\Lambda,\mu}$. Note that we did not introduce a different notation for the infinite vector $v$ when restricted to a finite $\Lambda$.
\qed
\end{proofsect}
%
%
%
 
\section*{Acknowledgments} Stefan Adams thanks David Brydges for many fruitful discussions on this subject. Quirin Vogel would like to thank Yacine Barhoumi for his advice on combinatorics.

\end{document}